\documentclass[preprint,12pt]{elsarticle}

\usepackage{lineno,hyperref,amsmath,amssymb,color}
\usepackage{graphics} 
\usepackage{epsfig}
\usepackage{graphicx}
\usepackage{epstopdf}
\modulolinenumbers[5]

\journal{Journal of \LaTeX\ Templates}
\usepackage{epstopdf}
\usepackage[colorinlistoftodos]{todonotes}
\usepackage{graphicx} 
\usepackage{psfrag} 
\usepackage{amsthm}
\usepackage{amsmath, amssymb, dsfont}
\usepackage{bm}
\usepackage{ulem}
\usepackage{geometry}
\newtheorem{theorem}{Theorem}[section]

\newtheorem{proposition}[theorem]{Proposition}

\newcommand\redout{\bgroup\markoverwith{\textcolor{red}{\rule[.5ex]{2pt}{0.8pt}}}\ULon}

\begin{document}
\begin{frontmatter}
\title{VMS spectral solution of two-dimensional advection-diffusion problems with anisotropic velocity}
%
\author{T. Chac\'on Rebollo, Soledad Fern\'andez-Garc\'ia\footnote{Corresponding author: soledad@us.es}}
\address{Dpto. EDAN\&IMUS, University of Seville, \\ Campus de Reina Mercedes, 41012 Sevilla, Spain}
\author{M. G\'omez-M\'armol}
\address{Dpto. EDAN, University of Seville, \\ Campus de Reina Mercedes, 41012 Sevilla, Spain}
%

%
%

\begin{abstract}
In this article, we extend the Variational Multi-scale method with spectral approximation of the sub-scales to two-dimensional advection-diffusion problems. The spectral VMS method is cast for low-order elements as a standard VMS method with specific stabilized coefficients, that are anisotropic in the sense that they depend on two grid P\'eclet numbers, each associated to a component of the advection velocity.
We compute the stabilized coefficients for grids of isosceles right triangles and right quadrilaterals, based upon the explicit computation of the eigen-pairs of the advection-diffusion operator with Dirichlet boundary conditions. To reduce the computing time, the stabilized coefficients are computed at the nodes of a grid in an off-line step, and then interpolated by a fast procedure in the on-line computation.

Finally, we present some numerical tests, first with constant velocity and after that with anisotropic variable velocity, in order to compare our results with those provided by other stabilization coefficients. We observe a relevant accuracy gain for moderately large grid P\'eclet numbers for variable advection velocity.
\end{abstract}


\end{frontmatter}

\section{Introduction} \label{sec:intro}
The Variational Multi-Scale method provides a general technique to deal with the instabilities generated by the Galerkin discretization of partial differential equations. These instabilities appear when the PDEs include terms of different derivation orders and the discretization parameters are not small enough. The discrete solutions present spurious oscillations when low-order operator terms are dominant. Such solutions are unreliable for technological and engineering applications (see Hughes (cf. \cite{Hughes_1995, HughesStewart_1995, HughesFMQ_1998}).  The pioneering stabilized method in the framework of finite element discretizations is the SUPG (Streamline Upwind Petrov-Galerkin) method (see \cite{Brooks_1982}). It consists in adding to the classical Galerkin formulation an extra term devoted to control the advection derivative. The SUPG method was followed by a large class of stabilized methods (Galerkin-Least Squares methods, adjoint --or unusual-- Galerkin-Least Squares methods, Orthogonal Sub-Scales (OSS) method, among others) all consisting in adding extra terms to the Galerkin formulation aiming to control one or several operator terms appearing in the equations. These methods where mainly applied to the numerical solution of incompressible and subsequently compressible flow equations, also proving that they provide a further stabilization of the discretization of the pressure gradient. An overview of most stabilized methods may be found in \cite{Hughes_1995}, while the OSS is introduced in\cite{codina12}.
\par
 The Variational Multi-Scale (VMS) formulation is based upon considering the effect on the small scales on the large ones at the discrete level. This effect is neglected by the Galerkin method. The small scales are driven by the residual of the large scales. A global stabilization effect is achieved (damping, at least partial, of spurious oscillations), due to a dissipative action of the small scales onto the large scales. To build a feasible VMS method, the small scale problem is further discretized by some kind of approximation. Within the VMS formulation, the above mentioned stabilized methods appear as ``partial" VMS methods that retain some of the diffusive terms that appear in the interaction large-small scales. The VMS methods have been successfully applied to many flow problems, and in particular to the building of  models of Large Eddy Simulation (LES) of turbulent flows, with remarkable accuracy (cf.\cite{hug112, john0612, Chaconlibro}). 
  \par
A possibility to approximately solve the small scales equation within the VMS framework is a local diagonalization of the PDE operator on each grid element. This leads to the Orthogonal Sub-Scales (OSS) method, introduced by Codina in \cite{codina12}. A further step within the diagonalization techniques is the use of spectral techniques, introduced in \cite{ChaconDia}. In this paper, the sub-grid scales are initially approximated by bubble functions on each grid element. The basic observation is that the eigen-pairs of the advection-diffusion operator may be calculated explicitly, which allows to analytically calculate the sub-grid scales by means of a spectral expansion on each grid element. A feasible VMS-spectral discretization is then built by truncation of this spectral expansion to a finite number of modes. For 1D advection-diffusion problems, the method with an odd number of modes satisfies the discrete maximum principle. It is found to be of 3rd. order with respect to the number of eigenmodes.
\par
In the present paper, we apply the spectral VMS method to 2D advection-diffusion problems. The spectral VMS method is re-written for low-order finite elements as a standard VMS method, with stabilized coefficients analytically computed from the eigen-pairs of the advection-diffusion operator. The novelty is that the stabilized coefficients depend on directional grid P\'eclet numbers, each associated to a coordinated direction. In this sense, we are dealing with an anisotropic VMS method. The eigen-pairs family of the advection-diffusion operator are explicitly constructed, for grids formed by either isosceles right triangles or straight quadrilaterals. To reduce the computing time, the stabilized coefficients are pre-computed at an off-line step at the nodes of an interpolation grid in the $Pe_{1h}$, $Pe_{2h}$ space, and then interpolated by a fast procedure in the on-line computation. This allows computing times quite close to the standard VMS methods.
\par
We perform some numerical tests to analyze the practical performance of the introduced method. For advection-diffusion problems with constant velocities the accuracy is quite close to that provided by the usual stabilized coefficients (either the extrapolation to 2D of the 1D coefficients, or the ones introduced by Codina in \cite{codina12}), for moderately large grid P\'eclet numbers (roughly below 10). For variable velocities, however, there is an appreciable gain that allows rates of error reduction in $L^2$ and $L^\infty$ norms up to $5$. This is particularly remarkable for realistic situations such as advection of passive scalars by a flow around a cylinder. For large grid P\'eclet numbers, the new stabilized coefficients fail as we are approximating the sub-grid scales by bubble functions, that vanish at the element boundaries. 
\par
In sum, the anisotropic spectral VMS method, for low-order elements, may be cast as a standard VMS method with specific stabilized coefficients, that take into account the anisotropy of the advection velocity. It provides gains of accuracy for advection-dominant flows at moderately large grid P\'eclet numbers, for variable advection velocity.
\par
The paper is outlined as follows. In Section \ref{sec:abstract}, we consider the abstract formulation of the problem. Section \ref{sec:2D} is devoted to applying the method to the two-dimensional advection-diffusion 
problem. The section is divided into two subsections, one where we consider a mesh of squares, and a second one where we consider a mesh of isosceles right triangles.  After that, in Section \ref{sec:stab} we compute the stabilization coefficients in an off-line phase. Section \ref{sec:numerical} is focused on numerical results to test the reliability of the method. Finally, Section \ref{sec:conclusions} is devoted to stating some conclusions and future work.
\section{Abstract Formulation} \label{sec:abstract}
We briefly remind in this section the abstract formulation of the Spectral VMS method for linear elliptic equations, which was developed in  \cite{ChaconDia}.

Let $X$ be a Hilbert space, $B$ a bilinear bounded and coercive form and $l\in X',$ the topological dual of $X.$ Consider the variational elliptic problem,
\begin{equation}\label{ecu_abstractDia}
\begin{array}{l}
B(U,V)=l(V), \quad \forall V\in X.
\end{array}
\end{equation}
Consider the decomposition,
$
X=X_h\oplus \tilde{X},
$
where $X_h$ is a sub-space of $X$ of finite dimension, and $\tilde{X}$ is a complementary, infinite-dimensional, sub-space of $X.$ 
Problem (\ref{ecu_abstractDia}) can be reformulated as
\begin{equation}\label{ecu_dec2}
\left\{\begin{array}{lr}
B(U_h+\tilde{U},V_h)=l(V_h),\quad\forall V_h\in X_h,&(a)\\ \noalign{\smallskip}
B(U_h+\tilde{U},\tilde{V})=l(\tilde{V}),\quad\forall  \tilde{V}\in \tilde{X},&(b)
\end{array}\right.
\end{equation}
where $$
\begin{array}{l}
U=U_h+\tilde{U}, \quad V=V_h+\tilde{V}, \quad \mbox{for } U_h, V_h\in X_h, \quad \tilde{U},\tilde{V}\in \tilde{X}.
\end{array}
$$
From equation (\ref{ecu_dec2})(b), it is possible to obtain the standard VMS reformulation of problem (\ref{ecu_abstractDia}), namely,
\begin{equation}\label{ecuVMS1}
\begin{array}{l}
B(U_h,V_h)+B(\Pi(R(U_h)),V_h)=l(V_h), \quad \forall V_h\in X_h,
\end{array}
\end{equation}
where $R$ is the residual of the large scales component and $\Pi$ is the static condensation operator. Note that  $\tilde{U}=\Pi(R(U_h)),$  (see  \cite{ChaconDia} for details).

Let us now consider that problem (\ref{ecu_abstractDia}) is the variational formulation of an elliptic PDE
$$
\begin{array}{l}
\mathcal{L}(U)=l,
\end{array}
$$
on a bounded domain $\Omega\subset\mathbf{R}^d$, that is, $\mathcal{L}:X\mapsto X'$ is the operator defined by
$$
\begin{array}{l}
\langle{\mathcal{L}}V,W\rangle=B(V,W),\quad \forall W\in X,
\end{array}
$$
and $X$ is a suitable Hilbert space of functions defined on $\Omega$. Our purpose is to do an spectral approximation of the small scales by using the standard VMS reformulation.

Given a triangulation $\mathcal{T}_h$ of the domain $\Omega,$ we can approximate the small scale space $\tilde{X}$ by 
\begin{equation*}\label{aprsmall}
\tilde{X}\simeq\bigoplus_{K\in\mathcal{T}_h}\tilde{X}_K,\quad\mbox{with } \tilde{X}_K=\{\tilde{V}\in\tilde{X}:\mathrm{supp}(\tilde{V})\subset K\}.
\end{equation*}
Hence,
$
\tilde{U}\simeq \sum_{K\in\mathcal{T}_h}\tilde{U}_K,\mbox{ with }\tilde{U}_K\in\tilde{X}_K,
$
and $\tilde{U}_K=\Pi_{K}(R(U_h)),$ where $\Pi_{K}$ denotes the restriction of operator $\Pi$ to  $\tilde{X}_K.$

Also, we denote  $\mathcal{L}_{K}$ the restriction of operator $\mathcal{L}$ to  $\tilde{X}_K.$ Now, we extend to dimension two the basic result of the spectral VMS method, studied in \cite{ChaconDia} in dimension one.
\begin{theorem}\label{maintheorem}
Let us assume that there exists a complete sub-set $\{\hat{z}_{js}^{(K)}\}_{j,s\in\mathbf{N}}$ on $\tilde{X}_K$ formed by eigenfunctions of the operator $\mathcal{L}_{K}$, which is an orthonormal system in $L^2_{p_{K}}(K)$ for some weight function $p_{K}\in C^1(\bar{K}).$ Then,
\begin{equation}\label{series1}
\begin{array}{l}
\tilde{U}_K=\displaystyle\sum_{j=1}^{\infty} \displaystyle\sum_{s=1}^{\infty}\beta_{js}^{(K)} \langle R(U_h),p_{K} \hat{z}_{js}^{(K)}\rangle \hat{z}_{js}^{(K)},\quad \mbox{with }\beta_{js}^{(K)} =(\lambda_{js}^{(K)})^{-1},\\
\end{array}
\end{equation}
where $\lambda_{js}^{(K)}$ is the eigenvalue of $\mathcal{L}_{K}$ associated to $\hat{z}_{js}^{(K)}.$
\end{theorem}
Finally, in order to obtain a feasible discretization, it is necessary to truncate series (\ref{series1}) to $M_1,M_2\geq 1$ addends and approximate problem (\ref{ecuVMS1}) by
\begin{equation}\label{ecuVMS_M}
\begin{array}{l}
B(U_{h,M_1,M_2},V_h)+B(\Pi^M(R(U_{h,M_1,M_2})),V_h)=l(V_h), \quad \forall V_h\in X_h,
\end{array}
\end{equation}
where the unknowns are $U_{h,M_1,M_2}\in X_h$ and the operator $\Pi^{M_1,M_2}$ is given by
$$
\Pi^{M_1,M_2}(\varphi)=\sum_{K\in\mathcal{T}_h}\Pi^{M_1,M_2}_{K}(\varphi), \,\mbox{with }\Pi^{M_1,M_2}_{K}(\varphi)=\sum_{j=1}^{M_1} \displaystyle\sum_{s=1}^{M_2}\beta_{js}^{(K)} \langle \varphi,p_{K} \hat{z}_{js}^{(K)}\rangle \hat{z}_{js}^{(K)}.
$$

Now that we have introduced the abstract formulation of the problem, let us focus our attention on the application to the two-dimensional advection-diffusion problem, which is the main goal of this article.

\section{Application to two-dimensional advection-diffusion problem}\label{sec:2D}
Consider the advection-diffusion problem in the unit square $\Omega=[0,1]\times[0,1]$,
\begin{equation}\label{EADS}
\left\{\begin{array}{l}
\mathbf{a}\cdot\nabla U-\mu \Delta U=f\quad \mbox{in }\Omega,\\  \noalign{\smallskip}
U=0 \quad \mbox{on }\partial \Omega,\\ 
\end{array}\right.
\end{equation}
where $\mathbf{a}\in L^\infty(\Omega)$ is a divergence-free given velocity field, $\mu>0$  is the diffusion coefficient and $f\in L^2(\Omega)$ is the source term.
The weak formulation of problem (\ref{EADS}) is given by
\begin{equation}\label{WEAD}
\left\{\begin{array}{l}
\mbox{Find }U\in H_0^1(\Omega) \quad\mbox{such that,}\\  \noalign{\smallskip}
(\mathbf{a}\cdot\nabla U,V)+\mu \,(\nabla U,\nabla V)=(f,V) \quad \forall V\in H_0^1(\Omega).\\  \noalign{\smallskip}
\end{array}\right.
\end{equation}
Now, we proceed to build the spectral VMS discretization (\ref{ecuVMS_M}) of problem (\ref{WEAD}).

Given a triangulation $\mathcal{T}_h$ of the domain $\Omega,$ let us assume that the velocity $\mathbf{a}$  is approximated in each sub-grid term by a constant value $\mathbf{a}_{K}$ on each element $K$.
Let us extend to dimension two a result from \cite{ChaconDia}, about the eigen-pairs of the advection-diffusion operator.
\begin{proposition}\label{propADR}
The couple $\left(\tilde{\omega}_{js}^{(K)},\lambda_{js}^{(K)}\right)$ is an eigenpair of the advection-diffusion 
operator $\mathcal{L}_{K}$ if and only if the couple $\left(\tilde{W}_{js}^{(K)},\sigma_{js}^{(K)}\right)$ is an eigenpair of the Laplace operator $-\Delta$ in $H_0^1(K),$ where
\begin{equation}\label{sigma1}
\begin{array}{l}
\tilde{\omega}_{js}^{(K)}=\psi^{(K)}\tilde{W}_{js}^{(K)}\quad\mbox{with}\quad \psi^{(K)}(\mathbf{x})=e^{\frac{1}{2\mu}(\mathbf{a}\cdot\mathbf{x})}\quad \\ \noalign{\medskip}
\mbox{and}\quad\lambda_{js}^{(K)}=\mu\left(\displaystyle\frac{|\mathbf{a}|^2}{4\mu^2}+\sigma_{js}^{(K)}\right),\forall {j,s}\in\mathbb{N}.
\end{array}
\end{equation}
\end{proposition}

The exact computation of the eigen-pairs of the Laplace operator can be done in the case of elements with simple geometrical forms, as it is the case of parallelepipeds and various type of triangles, for instance, equilateral \cite{Lame,McCartin,Pinsky} and isosceles right \cite{triangles}.  We focus our attention on two cases. In the first one, we consider a mesh of squares and in the second one, a mesh of isosceles right triangles.
\subsection{Mesh of squares}
Consider that the triangulation $\mathcal{T}_h$ of the domain $\Omega=[0,1]\times[0,1]$ is formed by $(N-1)\times(N-1)$ squares of side size $h=1/(N-1),N>2$ and denote a generic grid node as $(x_{l_1},y_{l_2})$ for $l_1,l_2=1,...,N$. The eigenvalue problem for the Laplace operator with Dirichlet boundary conditions in each element $K=(x_{l_1-1},x_{l_1})\times(y_{l_2-1},y_{l_2}),$  $l_1,l_2=2,...,N,$ is given by
\begin{equation}\label{SquareLaplace}
\left\{\begin{array}{l}
-\partial_{xx}W^{(K)}-\partial_{yy}W^{(K)}=\sigma W^{(K)} \quad\mbox{in }K, \\ \noalign{\smallskip}
W^{(K)}=0 \quad \mbox{on }\partial K.\\ 
\end{array}\right.
\end{equation}

Taking into account that in the 1D problem with uniformly spaced nodes with $x_l-x_{l-1}=h,$  the eigen-pairs are
\begin{equation*}\label{sigma2}
\tilde{W}_j^{(K)}(x)=\sin\left(\sqrt{\sigma_j^{(K)}}(x-x_l)\right)\mbox{ with } \sigma_j^{(K)}=\left(\frac{j\pi}{h}\right)^2 \mbox{ for }  l=1,...,N \mbox{ and } j\in\mathbb{N},
\end{equation*}
one can see that the solutions of problem (\ref{SquareLaplace}) are given by
\begin{equation}\label{sigma2D}
\begin{array}{l}
W_{js}^{(K)}(x,y)=\sin\left(\sqrt{\sigma_j^{(K)}}(x-x_{l_1})\right) \sin\left(\sqrt{\sigma_s^{(K)}}(y-y_{l_2})\right), \\   \noalign{\smallskip}
\mbox{ with } \sigma_{j}^{(K)}=\left(\displaystyle\frac{j\pi}{h}\right)^2, \sigma_{s}^{(K)}=\left(\displaystyle\frac{s\pi}{h}\right)^2 \mbox{ for }  l_1,l_2=1,...,N \mbox{ and } j,s \in\mathbb{N},
\end{array}
\end{equation}
where the corresponding eigenvalue is $\sigma_{js}^{(K)}=\sigma_{j}^{(K)}+\sigma_{s}^{(K)}.$ Hence, from Proposition \ref{propADR}, the eigen-pairs of the 2D advection-diffusion operator with Dirichlet boundary conditions on $K$ are given by
\begin{equation}\label{sigma3}
\begin{array}{l}
\tilde{\omega}_{js}^{(K)}(x,y)=e^{\frac{1}{2\mu}(\mathbf{a}\cdot\mathbf{x})} \sin\left(\sqrt{\sigma_j^{(K)}}(x-x_{l_1})\right) \sin\left(\sqrt{\sigma_s^{(K)}}(y-y_{l_2})\right) \\ \noalign{\medskip}
\mbox{and}\quad\lambda_{js}^{(K)}=\mu\left(\displaystyle\frac{|\mathbf{a}|^2}{4\mu^2}+\left(\displaystyle\frac{j\pi}{h}\right)^2+\left(\displaystyle\frac{s\pi}{h}\right)^2\right),\mbox{ for }  l_1,l_2=1,...,N \mbox{ and } j,s \in\mathbb{N}.
\end{array}
\end{equation}

The following result is an extension to two dimensions to that proven in \cite{ChaconDia} in dimension one. It shows that Theorem \ref{maintheorem} is satisfied in our setting.
\begin{theorem}\label{ortonormal2}
The sequence $\{\tilde{z}_{js}^{(K)}\}_{j,s \in\mathbf{N}}$ is complete in $H_0^1(K)$ and orthonormal in $L^2_{p_K}(K),$ where
$$\tilde{z}_{js}^{(K)}(x,y)=\tilde{\omega}_{js}^{(K)}/\|\tilde{\omega}_{js}^{(K)}\|_{p_K}= \frac{2 e^{\frac{1}{2\mu}(\mathbf{a}\cdot\mathbf{x})} }{h}\sin\left(\sqrt{\sigma_j^{(K)}}(x-x_{l_1})\right) \sin\left(\sqrt{\sigma_s^{(K)}}(y-y_{l_2})\right) $$
for $  l_1,l_2=1,...,N$ and $j,s \in\mathbb{N}$ and the weight functions are given by
\begin{equation}\label{pesos}
p_{K}= (\psi^{(K)})^{-2}=e^{-\frac{1}{\mu}(\mathbf{a}\cdot\mathbf{x})},
\end{equation}
where $\tilde{\omega}_{js}^{(K)}$ and $\psi^{(K)}$ are provided in expression (\ref{sigma1}).
\end{theorem}
Thus, from Theorem \ref{maintheorem} and expression (\ref{sigma3}), the stabilization coefficients in the case of a triangulation of the domain formed by squares of side size $h$ are given by
\begin{equation}
\beta_{js}^{(K)}=\frac{h^2}{\mu\left(Pe_h^2+\pi^2(j^2+s^2) \right)},\mbox{ for } j,s \in\mathbb{N},
\end{equation}
being $Pe_h=\frac{h|\mathbf{a}|}{2\mu}$ the element P\'eclet number.

\subsection{Mesh of isosceles right triangles}
\begin{figure}[h!]
\begin{center}
\includegraphics[width=0.5\linewidth]{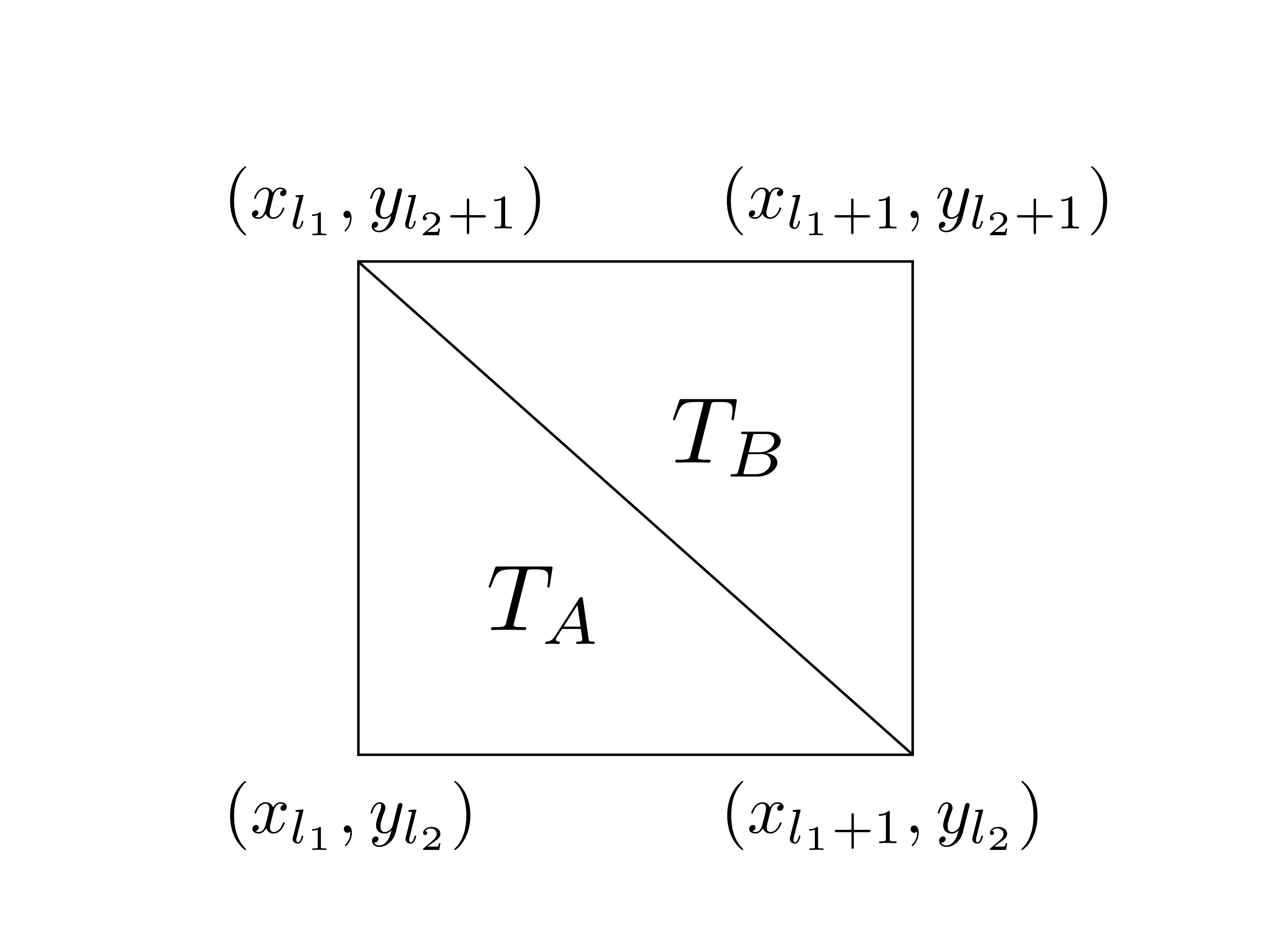}
\end{center}
\caption{\label{fig.trig} The two types of isosceles right triangles considered.}
\end{figure}

Consider that the triangulation $\mathcal{T}_h$ of the domain $\Omega=[0,1]\times[0,1]$ is formed by isosceles right triangles splitting squares of side size $ h=1/(N-1),N>2 $ by two and denote a generic grid node as $(x_{l_1},y_{l_2})$ for $l_1,l_2=1,...,N$. For the sake of clarity, we denote $T_A$ the triangle below the diagonal and $T_B$ the triangle above it, see Fig. \ref{fig.trig}.

In order to obtain the eigen-pairs of the advection-diffusion problem with Dirichlet boundary conditions in a generic triangle $T$ of the triangulation, we begin by computing them in the reference triangle 
\begin{equation}\label{triangref}
\hat{T}=\{(\hat{x},\hat{y}):\hat{x}\in(0,1),y\in(0,1-\hat{x})\}.
\end{equation}
As a previous step, we need to know the eigen-pairs of the Laplace operator with Dirichlet boundary conditions in  $\hat{T}.$
From \cite{triangles} we know that these eigenfunctions are given by
\begin{equation}\label{sigma2D}
\begin{array}{l}
\hat{W}_{js}(\hat{x},\hat{y})=\sin\left(j\pi \hat{x}\right) \sin\left(s\pi \hat{y}\right)-\sin\left(s\pi(1-\hat{x})\right) \sin\left(j\pi(1-\hat{y})\right), \\   \noalign{\smallskip}
\mbox{ with the eigenvalues  } \hat{\sigma}_{js}=\pi^2(j^2 +s^2)\mbox{ for }  j,s \in\mathbb{N}.
\end{array}
\end{equation}
Hence, from Proposition \ref{propADR}, the eigen-pairs of the 2D advection-diffusion operator with Dirichlet boundary conditions on the reference triangle $\hat{T}$ are given by
\begin{equation*}\label{sigma4}
\begin{array}{l}
\hat{\omega}_{js}(\hat{x},\hat{y})=e^{\frac{1}{2\mu}(\mathbf{a}\cdot\hat{\mathbf{x}})}\left(\sin\left(j\pi \hat{x}\right) \sin\left(s\pi \hat{y}\right)-\sin\left(s\pi(1-\hat{x})\right) \sin\left(j\pi(1-\hat{y})\right)\right) \\ \noalign{\medskip}
\mbox{and}\quad\hat{\lambda}_{js}=\mu\left( \displaystyle\frac{|\mathbf{a}|^2}{4\mu^2}+\pi^2(j^2+s^2)\right),\mbox{ for }  j,s \in\mathbb{N}.
\end{array}
\end{equation*}

The following result shows that Theorem \ref{maintheorem} is satisfied in our setting.
\begin{theorem}\label{ortonormal}
The sequence $\{\hat{z}_{js}\}_{j,s \in\mathbf{N}}$ is complete in $H_0^1(\hat{T})$ and orthonormal in $L^2_{p_{\hat{T}}}(\hat{T}),$ where
\begin{equation}
\begin{array}{rcl}
\hat{z}_{js}(\hat{x},\hat{y})&=&\hat{\omega}_{js}/\|\hat{\omega}_{js}\|_{p_K} \\ \noalign{\medskip}
&=&2 e^{\frac{1}{2\mu}(\mathbf{a}\cdot\hat{\mathbf{x}})}\hat{W}_{js}(\hat{x},\hat{y}) \mbox{ for }  j,s \in\mathbb{N},
\end{array}
\end{equation}
being the weight functions those provided in expression (\ref{pesos}), with $\hat{W}_{js}$ is given in expression (\ref{sigma2D}) and $\hat{T}$ the reference triangle given in (\ref{triangref}).
\end{theorem}
Once we have found the eigen-pairs of the 2D advection-diffusion operator with Dirichlet boundary conditions in the reference triangle $\hat{T},$ we proceed to compute 
them in a generic triangle $T\in \mathcal{T}_h$, by considering the appropriate change of variables. 

Remember that the triangulation $\mathcal{T}_h$ of the domain $\Omega=[0,1]\times[0,1]$ is formed by isosceles right triangles splitting squares of side size $ h $ by two and that we denote a generic grid node as $(x_{l_1},y_{l_2})$ for $l_1,l_2=1,...,N$. As we have commented before, we distinguish two different type of triangles in the mesh, depending if they are below the diagonal or above it. In particular, in a generic square $K=(x_{l_1-1},x_{l_1})\times(y_{l_2-1},y_{l_2}),$  $l_1,l_2=2,...,N,$ we consider the two triangles
\begin{equation}
T_A=\{(x,y):x\in(x_{l_1-1},x_{l_1}), y\in (y_{l_2-1},x_{l_1}-x+y_{l_2})\}
\end{equation}
and
\begin{equation}
T_B=\{(x,y):x\in(x_{l_1-1},x_{l_1}), y\in (x_{l_1}-x+y_{l_2},y_{l_2})\}
\end{equation}
being $T_A$ the triangle below the diagonal and $T_B$ the triangle above it, see Fig. \ref{fig.trig} again.
Consider the applications from triangles $T_A$ and $T_B,$ respectively, to the reference element, that is,
\begin{equation*}\label{GA}
\begin{array}{cccc}
G_A:&T_A&\longrightarrow& \hat{T}\\ 
&(x,y)&\longmapsto&(\hat{x}_A,\hat{y}_A)=\displaystyle\left(\frac{x-x_{l_1-1}}{h},\frac{y-y_{l_2-1}}{h}\right)
\end{array}
\end{equation*}
and
\begin{equation*}\label{GB}
\begin{array}{cccc}
G_B:&T_B&\longrightarrow& \hat{T}\\ 
&(x,y)&\longmapsto&(\hat{x}_B,\hat{y}_B)=\displaystyle\left(\frac{x_{l_1}-x}{h},\frac{y_{l_2}-y}{h}\right).
\end{array}
\end{equation*}
Now, it is easy to see that the following proposition holds.
\begin{proposition}
The eigenfunctions of the two-dimensional advection-diffusion operator with Dirichlet boundary conditions on triangles $T_A$ and $T_B$ are given by
\begin{equation}\label{eigTA}
\begin{array}{l}
\omega_{js}^A(x,y)=\displaystyle\exp\left(\frac{a_1}{2\mu}\frac{x-x_{l_1-1}}{h}+\frac{a_2}{2\mu}\frac{y-y_{l_2-1}}{h}\right)f^A_{js}(x,y), \mbox{ for }  l_1,l_2=1,...,N \mbox{ and } j,s \in\mathbb{N}\\ \noalign{\medskip}
\end{array}
\end{equation}
and
\begin{equation}\label{eigTB}
\begin{array}{l}
\omega_{js}^B(x,y)=\displaystyle\exp\left(\frac{a_1}{2\mu}\frac{x-x_{l_1-1}}{h}+\frac{a_2}{2\mu}\frac{y-y_{l_2-1}}{h}\right)f^B_{js}(x,y), \mbox{ for }  l_1,l_2=1,...,N \mbox{ and } j,s \in\mathbb{N},\\ \noalign{\medskip}
\end{array}
\end{equation}
where
\begin{equation*}\label{fA}
\begin{array}{l}
f^A_{js}(x,y)=\displaystyle\sin\left(j\pi\frac{x-x_{l_1-1}}{h}\right) \sin\left(s\pi \frac{y-y_{l_2-1}}{h}\right)-\sin\left(s\pi \frac{x_{l_1}-x}{h}\right) \sin\left(j\pi\frac{y_{l_2}-y}{h}\right)
\end{array}
\end{equation*}
and
\begin{equation*}\label{fB}
\begin{array}{l}
f^B_{js}(x,y)=\displaystyle \sin\left(j\pi \frac{x_{l_1}-x}{h}\right) \sin\left(s\pi\frac{y_{l_2}-y}{h}\right)-\sin\left(j\pi\frac{x-x_{l_1-1}}{h}\right) \sin\left(s\pi \frac{y-y_{l_2-1}}{h}\right),
\end{array}
\end{equation*}
respectively, and the eigenvalues are given by
\begin{equation}\label{eigvTA}
\begin{array}{l}
\lambda_{js}^{T}=\mu\left(\displaystyle\frac{|\mathbf{a}|^2}{4\mu^2}+\left(\displaystyle\frac{j\pi}{h}\right)^2+\left(\displaystyle\frac{s\pi}{h}\right)^2\right),\mbox{ for }  j,s \in\mathbb{N}.
\end{array}
\end{equation}
\end{proposition}
Hence, the following result holds,
\begin{theorem}\label{ortonormal}
The sequences $\{z_{js}^{A}\}_{j,s \in\mathbf{N}}$ and $\{z_{js}^{B}\}_{j,s \in\mathbf{N}}$ are complete in $H_0^1(T^A)$ and $H_0^1(T^B),$ respectively, and orthonormal in $L^2_{p_A}(T^A)$ and
$L^2_{p_B}(T^B),$ respectively, where
$$z_{js}^{A}(x,y)=\omega_{js}^{A}/\|\omega_{js}^{A}\|_{p_A}= \frac{2}{h} \omega_{js}^A(x,y)\quad\mbox{and}\quad z_{js}^{B}(x,y)=\omega_{js}^{B}/\|\omega_{js}^{B}\|_{p_B}= \frac{2}{h} \omega_{js}^B(x,y),$$
for $j,s \in\mathbb{N},$ where $\omega_{js}^{A}, \omega_{js}^{B}$ are given in expressions (\ref{eigTA}) and (\ref{eigTB}) and the weight functions are given by
\begin{equation*}\label{pesosAB}
p_{A}=e^{-\frac{1}{\mu}(\mathbf{a}\cdot\hat{\mathbf{x}}_A)}\quad\mbox{and}\quad p_{B}= e^{-\frac{1}{\mu}(\mathbf{a}\cdot\hat{\mathbf{x}}_B)}.
\end{equation*}
\end{theorem}
Thus, from Theorem \ref{maintheorem} and expression (\ref{eigvTA}) the stabilization coefficients in the case of a triangulation of the domain formed by isosceles right triangles of side size $h$ are given by
\begin{equation}\label{beta}
\beta_{js}^{(K)}=\frac{h^2}{\mu\left(Pe_h^2+\pi^2(j^2+s^2) \right)},\mbox{ for } j,s \in\mathbb{N},
\end{equation}
being $Pe_h=\frac{h|\mathbf{a}|}{2\mu}$ the element P\'eclet number.

Now that we have computed the eigen-pairs of the 2D advection-diffusion operator with Dirichlet boundary conditions in two different type of meshes, we can explicitly compute the stabilization coefficients of the VMS-spectral method, in both cases.
 
\section{Computation of the stabilization coefficients}\label{sec:stab}
Here, we extend to dimension two the relationship between the VMS-spectral method (\ref{ecuVMS_M}) and the usual VMS methods, which has been studied in 1D in \cite{ChaconDia}. Note that, due to element-wise regularity of $U_{h,M_1,M_2}$ and $V_h,$ we know that
\begin{equation*} \label{metsvms}
\begin{array}{rcl}
B(\Pi^{M_1,M_2}(R(U_{h,M_1,M_2})),V_h)&=&\displaystyle  \sum_{K\in\mathcal{T}_h}\sum_{j=1}^{M_1} \sum_{s=1}^{M_2}\beta_{js}^{(K)}\langle l- \mathcal{L} U_{h,M_1,M_2},p_{K} \hat{z}_{js}^{(K)} \rangle B(\hat{z}_{js}^{(K)},V_h)\\ \noalign{\smallskip}
&=& \displaystyle  \sum_{K\in\mathcal{T}_h}\sum_{j=1}^{M_1}  \sum_{s=1}^{M_2}\beta_{js}^{(K)}( l- \mathcal{L} U_{h,M_1,M_2},p_{K} \hat{z}_{js}^{(K)} )_K (\mathcal{L}^*V_h,\hat{z}_{js}^{(K)})_K.
\end{array}
\end{equation*}
If we consider $Q^1$ elements in the case of a mesh of squares, $P^1$ elements in the case of triangles, and we take into account that the velocity $\mathbf{a}$  is approximated in the sub-grid term by a constant value $\mathbf{a}_{K}$ on each element $K$, it follows that
\begin{equation}
\begin{array}{rcl}
B(\Pi^{M_1,M_2}(R(U_{h,M_1,M_2})),V_h)&=& B^S(U_{h,M_1,M_2},V_h)-l^S(V_h), 
\end{array}
\end{equation}
where
\begin{equation}
B^S(U_{h,M_1,M_2},V_h)=\displaystyle \sum_{K\in\mathcal{T}_h}\sum_{j=1}^{M_1}  \sum_{s=1}^{M_2} \beta_{js}^{(K)} (\mathbf{a}_K\cdot \nabla U_{h,M_1,M_2},p_{K} \hat{z}_{js}^{(K)})_K(\mathbf{a}_K\cdot \nabla V_h,\hat{z}_{js}^{(K)})_K
\end{equation}
and 
 $$
l^S(V_h)=\displaystyle \sum_{K\in\mathcal{T}_h} \sum_{j=1}^{M_1}  \sum_{s=1}^{M_2}\beta_{js}^{(K)} (l,p_{K} \hat{z}_{js}^{(K)})_K(\mathbf{a}_K\cdot \nabla V_h,\hat{z}_{js}^{(K)})_K.
$$
Thus, we can rewrite
\begin{equation}\label{equality}
B^S(U_{h,M_1,M_2},V_h)\simeq\displaystyle \sum_{K\in\mathcal{T}_h} \tau^{(K)}_{M_1,M_2}(\mathbf{a}_K\cdot \nabla U_{h,M_1,M_2},\mathbf{a}_K\cdot \nabla V_h)_K
\end{equation}
and 
$$
l^S(V_h)=\displaystyle \sum_{K\in\mathcal{T}_h} \tau^{(K)}_{M_1,M_2}(l,\mathbf{a}_K\cdot \nabla V_h)_K,
$$ 
where
\begin{equation}\label{tau_stac1}
 \tau^{(K)}_{M_1,M_2} =\displaystyle \frac{1}{|K|} \sum_{j=1}^{M_1}  \sum_{s=1}^{M_2} \beta_{js}^{(K)}  \left(\int_K  p_{K} \hat{z}_{js}^{(K)} \right)  \left(\int_K  \hat{z}_{js}^{(K)} \right),
\end{equation}
being the approximation (\ref{equality}) exact in the case of $P^1$ elements. 
Thus, the spectral method can be cast as a SUPG method with the stabilized coefficients $\tau^{(K)}_{M_1,M_2}$ and the standard stability and error analysis applies.

Finally, let us effectively compute the stabilization coefficients given in expression (\ref{tau_stac1}) in the case of a grid of isosceles right  triangles.
Assuming that $|K|=h^2/2,$ and taking into account the expression of $\beta_{js}^{(K)},$ given in expression (\ref{beta}), it follows that 
\begin{equation}
 \tau^{(K)}_{M_1,M_2}  =\displaystyle \frac{2}{h^2}\sum_{j=1}^{M_1}  \sum_{s=1}^{M_2}\left(\frac{\mu}{h^2}\left(\pi^2(j^2+s^2)+Pe_h^2\right)\right)^{-1} \left(\int_K  p_{K} \hat{z}_{js}^{(K)} \right)  \left(\int_K  \hat{z}_{js}^{(K)} \right),
\end{equation}
where we recall that $Pe_h=\frac{h |\mathbf{a}|}{2\mu}$ is the grid P\'eclet number. Now, a change of variables to the reference element let us conclude that
\begin{equation}\label{tauP}
 \tau^{(K)}_{M_1,M_2}  =\displaystyle \frac{8h^2}{\mu}\,\psi_{M_1,M_2}({Pe_h}_1,{Pe_h}_2),
\end{equation}
where 
$$
\psi_{M_1,M_2}({Pe_h}_1,{Pe_h}_2)=\sum_{j=1}^{M_1}  \sum_{s=1}^{M_2}\left(\pi^2(j^2+s^2)+Pe_h^2\right)^{-1} \hat{I}_{j,s}^-({Pe_h}_1,{Pe_h}_2)  \hat{I}_{j,s}^+ ({Pe_h}_1,{Pe_h}_2),
$$
being 
$$
 \hat{I}_{j,s}^-({Pe_h}_1,{Pe_h}_2)=\int_{\hat{T}}e^{-{Pe_h}_1(\hat{x}-1)-{Pe_h}_2(\hat{y}-1)}\hat{W}_{js}(\hat{x},\hat{y})d\hat{x}d\hat{y}
$$
and
$$
 \hat{I}_{j,s}^+({Pe_h}_1,{Pe_h}_2)=\int_{\hat{T}}e^{{Pe_h}_1(\hat{x}-1)+{Pe_h}_2(\hat{y}-1)}\hat{W}_{js}(\hat{x},\hat{y})d\hat{x}d\hat{y},
$$
with $\hat{W}_{js}$ given in expression (\ref{sigma2D}), ${Pe_h}_1=\frac{h a_1}{2\mu}$ and ${Pe_h}_2=\frac{h a_2}{2\mu},$ the grid P\'eclet numbers in each direction.

Thus, from expression (\ref{tauP}), it follows that it is possible to calculate the stabilization coefficients at suitable interpolation nodes of the plane ${Pe_h}_1,{Pe_h}_2$ in the off-line phase once for all. Then, these are interpolated by a fast procedure in the on-line phase, i.e., when actually computing the solution of the VMS-spectral method.

To measure the quality of the computed spectral stabilization coefficient, we will make use of the fact that the stabilization coefficients that yield the optimal diffusion in the 1D advection-diffusion problem are well-known and are given by
\begin{equation}\label{tau1D}
 \tau_{1D} =\displaystyle \frac{  \mu}{a_K^2}\,\varphi(Pe_h),\,\mbox{ with  } \varphi(P)=(P \coth(P)-1),
\end{equation}
being $a_K$ the one dimensional velocity in the $K$ element.
Thus, we are going to compare the stabilization coefficient obtained from formula (\ref{tauP}) in 2D advection-diffusion problems with zero velocity in one of the directions (that is, either ${Pe_h}_1=0$ or  ${Pe_h}_2=0$), with the optimal stabilization coefficients in 1D.  Assume, for instance that ${Pe_h}_2=0.$ Hence, $\mathbf{a}_K=(a_{1_K},0).$ Then, if we numerate the grid elements on the triangulation of the unit square in rows and columns and we denote as $K_{m,n}$ the element in row $m$ and column $n,$ from expression (\ref{equality}), we can write 
\begin{equation}\label{equality2}
\begin{array}{rcl}
B^S(U_{h,M_1,M_2},V_h)&\simeq&\displaystyle \sum_{m} \tau_{2D} \sum_{n} ((a_{1_K},0)\cdot \nabla U_{h,M_1,M_2},(a_{1_K},0)\cdot \nabla V_h)_{K_{m,n}} \\ \noalign{\smallskip}
&=&\tau_{2D}a_{1_K}^2\displaystyle   \sum_{m} (\partial_x U_{h,M_1,M_2}, \partial_x V_h)_{[x_n,x_{n+1}]}.
\end{array}
\end{equation}
Thus, it makes sense to compare the optimal numerical diffusion in 1D, which is given by $\nu_{num}^{1D}=a^2 \tau_{1D},$ with the obtained numerical diffusion from the 2D problem in expression (\ref{equality2}), that is, $\nu_{num}^{2D}=\tau_{2D}a^2,$
where we have chosen $a_K=a_{1_K}=a.$ That is, we have to compare $\tau_{1D}$ given in expression (\ref{tau1D}) with $\tau_{2D}$ given in (\ref{tauP}) when   ${Pe_h}_1=P$ and ${Pe_h}_2=0.$ This is equivalent to compare
$(\varphi(P))/(4P)$ with $8P\psi_{M_1,M_2}(P,0).$ Analogously, if we consider a case with  ${Pe_h}_1=0$ and ${Pe_h}_2=P,$ we obtain that we have to compare with $8P\psi_{M_1,M_2}(0,P).$

\begin{figure}[h!]
\begin{center}
\includegraphics[width=0.5\linewidth]{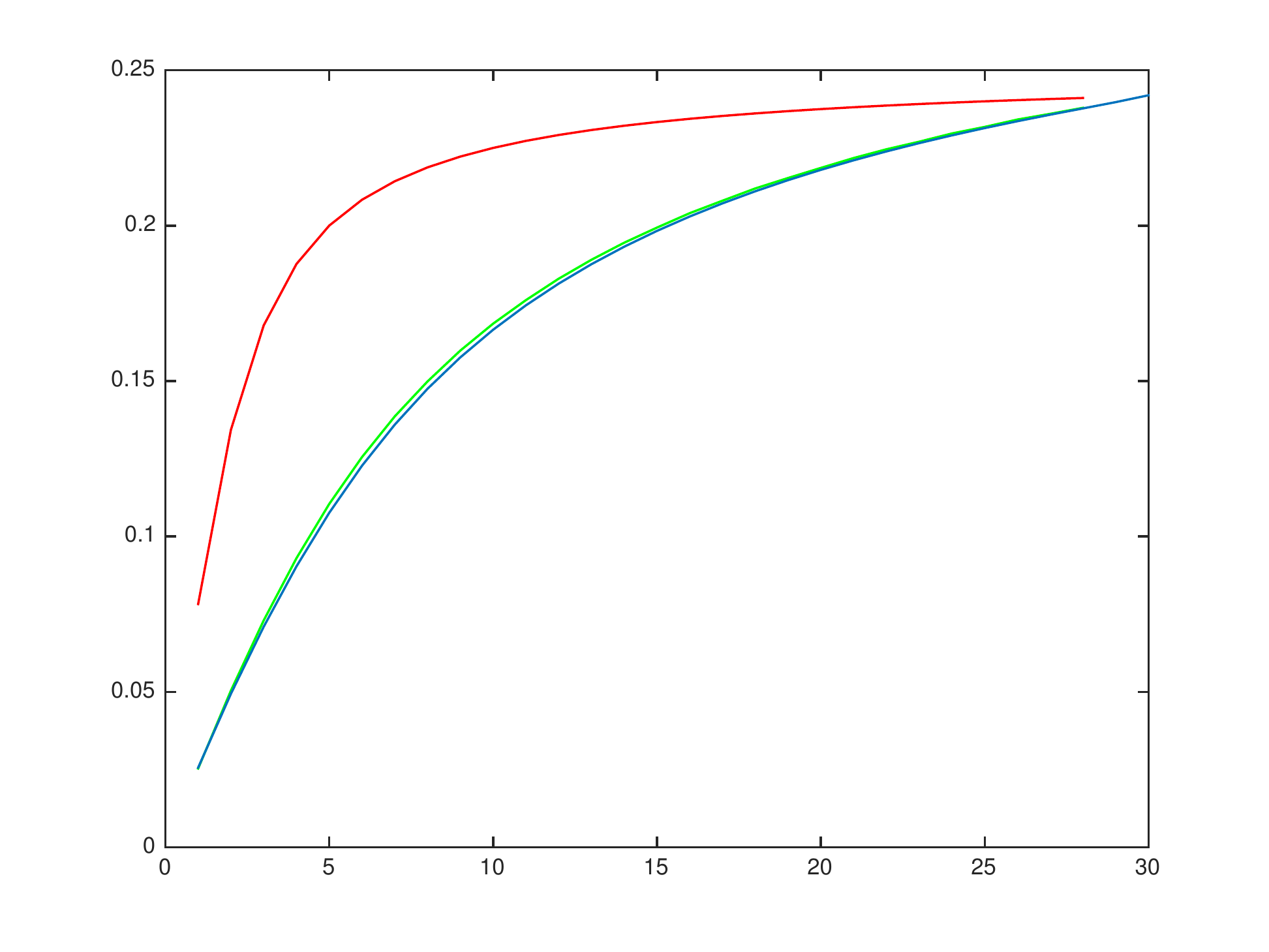}
\end{center}
\caption{\label{fig_stab1D} Comparison between stabilized 1D and 2D coefficients. Function $(\varphi(P))/(4P)$ in red, $8P\psi_{M_1,M_2}(P,0)$ in green and $8P\psi_{M_1,M_2}(0,P)$ in blue, with $P\in(0,30)$ and $M_1,M_2$ large enough.}
\end{figure}

In Fig. \ref{fig_stab1D}, we have represented function $(\varphi(P))/(4P)$ in red, $8P\psi_{M_1,M_2}(P,0)$ in green and $8P\psi_{M_1,M_2}(0,P)$ in blue, with $P\in(0,30)$ and $M_1,M_2$ large enough. 

We can observe a difference that decreases as $P$ increases between the stabilization coefficients obtained in 2D and in 1D. We know that the stabilization coefficient in 1D is optimal, due to the fact that, in this case, the decomposition of the Hilbert space into the large scales space plus the small scales space is exact. Nevertheless, this is not the case in the 2D framework, where we would need to add subscales on the element boundaries in order to have an exact decomposition of the 2D space. 
Moreover, as it can be observed, while function $\varphi(P)/(4P)$ has an asymptote when $P$ grows, the functions coming from the 2D expression, $8P\psi_{M_1,M_2}(P,0)$ and $8P\psi_{M_1,M_2}(0,P),$ do not behave in the same way. Therefore, we consider that the obtained 2D coefficients can be used for problems with ``moderate'' P\'eclet numbers.

\section{Numerical Results}\label{sec:numerical}
In this section, we are going to compare the results provided by using the spectral stabilized coefficients, to those obtained through other methods. 
In particular, we will consider two other coefficients. 

First, 
we consider the generalization to 2D of the optimal stabilization coefficient in 1D \cite{christie,johnkno}, that is,
\begin{equation}\label{tauL}
 \tau_{G}^{(K)} =\displaystyle \frac{  \mu}{\|\mathbf{a}_K\|^2}\,\varphi(Pe_h),\,\mbox{ with  } \varphi(P)=(P \coth(P)-1),
\end{equation}
being $Pe_h=\frac{h_K|\mathbf{a}|}{2\mu}$ the element P\'eclet number.

Second, the stabilization coefficient through orthogonal sub-scales in finite element methods proposed by Codina in \cite{codina12}, namely,
\begin{equation}\label{taucodina}
\tau^{(K)}_C=\left(\left(4\frac{\mu}{h_K^2}\right)^2+\left(2\frac{|\mathbf{a}_K|}{h_K}\right)^2\right)^{-1/2}.
\end{equation}
%

We divide the results into three subsections. In the first one, we focus on tests with constant velocity. In the second one, we consider a test with anisotropic velocity in a regular mesh. Finally, we consider a flow around a cylinder velocity and an irregular mesh. 

\subsection{Tests with constant velocity}
In this section, we will denote
$\tau_A$ the stabilization coefficient of the VMS-spectral method given in (\ref{tauP}), $\tau_G$ the stabilization coefficient given in (\ref{tauL}) and $\tau_C$ the stabilization coefficient proposed by Codina given in (\ref{taucodina}).

As in the computation of the VMS-spectral method, differently to previous methods, the direction of the velocity field is taken into account, we proceed here to compare the performance of the VMS-spectral method in comparison with stabilization coefficients $\tau_G$ and $\tau_C$ when we consider velocities in different directions. To to that, we
consider the advection-diffusion problem (\ref{EADS}) in the unit square with Dirichlet boundary conditions, with diffusion coefficient $\mu=1$ and source term $f(x,y)=\sin(\pi x)\cos(\pi y)$. We take a triangular mesh of isosceles right triangles of side size $h=1/80$ and consider constant velocities of the form  $\mathbf{a}= (800\sqrt{2}\cos{\alpha},800\sqrt{2}\sin{\alpha}),$ being $\alpha=n\pi/10,$ for $n=0,2,\ldots,18.$ In these cases, 
the global P\'eclet number is ${Pe_h}=7.07,$ while the directional P\'eclet numbers ${Pe_h}_1=(a_1 h)/(2\mu),$ ${Pe_h}_2=(a_2 h)/(2\mu)$ vary from $-7,07$ to $7.07.$

\begin{figure}[h!]
\begin{center}
\begin{tabular}{ll}(a)  Errors in norm $L^2$ &(b)  Errors in norm $L^\infty$\\
\includegraphics[width=0.5\linewidth]{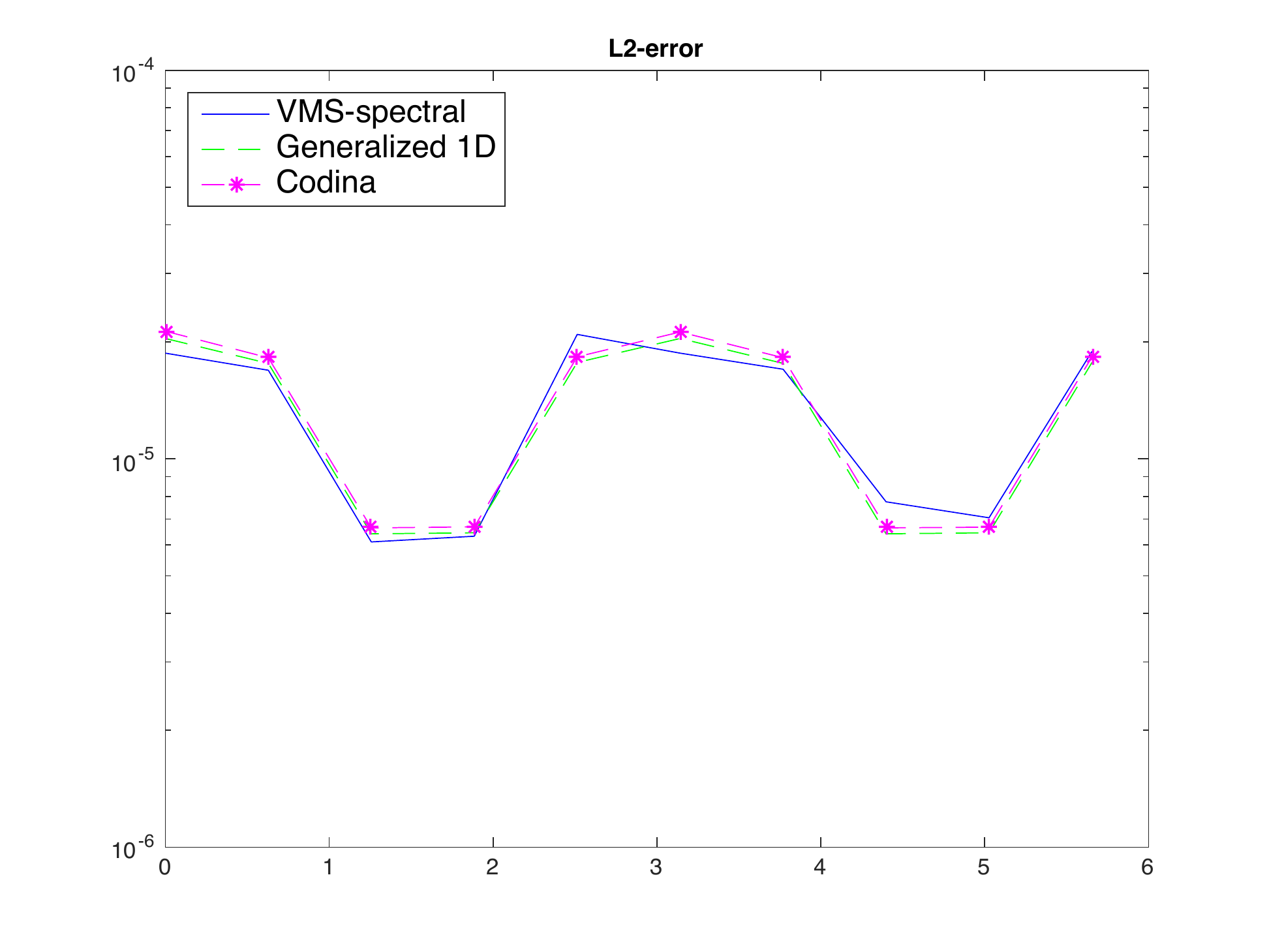}&\includegraphics[width=0.5\linewidth]{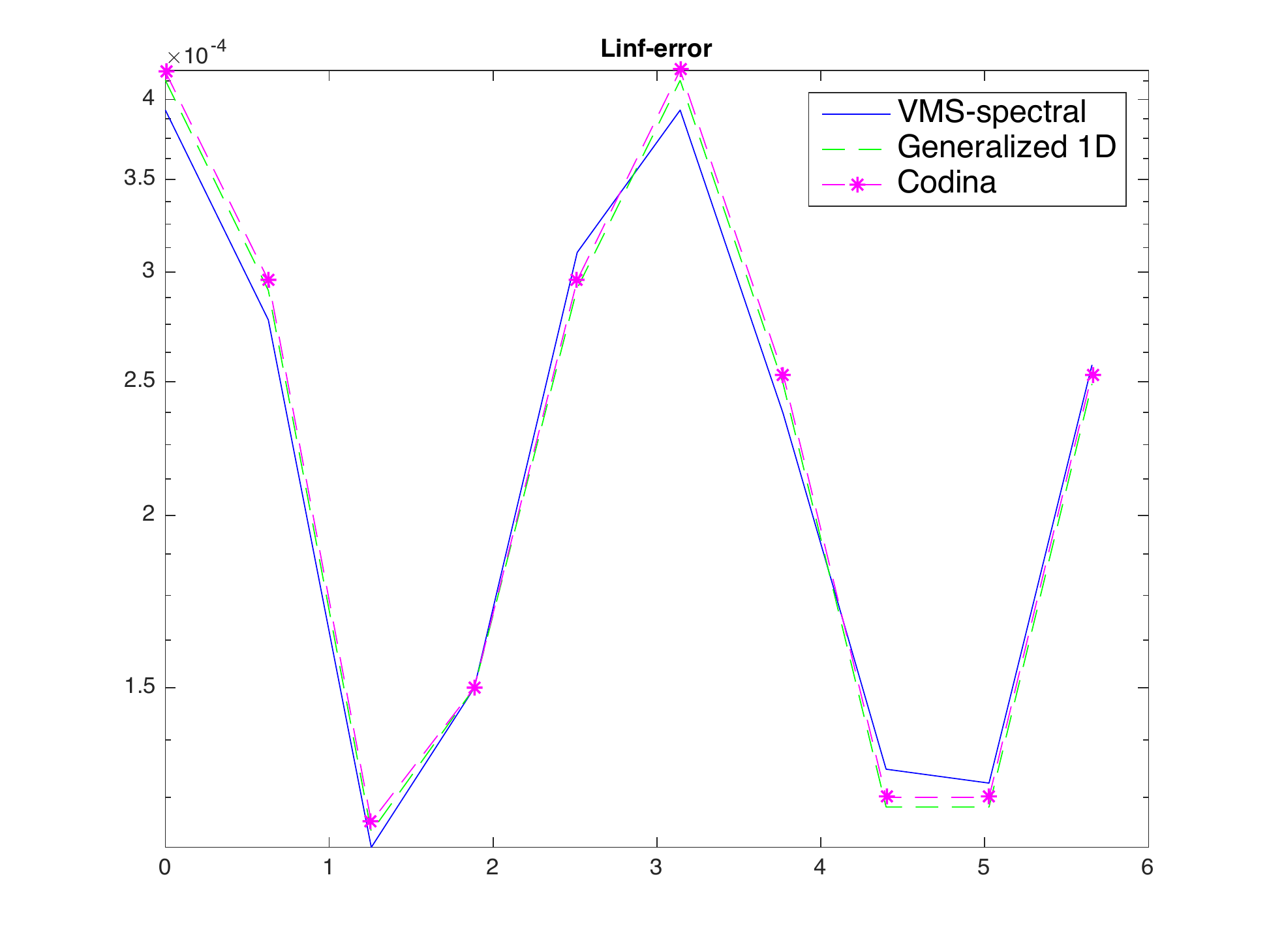}
\end{tabular}
\end{center}
\caption{\label{errores} Representation of the errors obtained with different schemes in norms $L^2$ (panel (a)) and $L^\infty$, (panel (b)) for $\alpha=n\pi/10,$ with $n=0,2,\ldots,18.$}
\end{figure}

In Figure \ref{errores}, we represent for each $n=0,2,\ldots,18,$
the error in $L^2$ and $L^\infty$ norms of the stabilized solution through the VMS-spectral stabilized coefficients (with subindex $A$), through the generalized 1D stabilized coefficients (with subindex $G$), and through Codina stabilized coefficients (with subindex $C$), by comparing the stabilized solution with the solution of the problem in a grid with 6561 degrees of freedom. These errors roughly behave as periodic functions with period $\pi$ in the three cases. All three stabilized coefficients settings provide quite similar error levels, no one of them appears to provide a better accuracy in all cases.

In subsequent figures, we consider the cases $n=4$ ($\alpha=2\pi/5$) and $n=10$ ($\alpha=\pi$) that respectively correspond to the smallest and largest errors. 

\begin{figure}[h!]
\begin{center}
\begin{tabular}{ll}(a)   VMS-spectral method&(b)  Generalized 1D method\\
\includegraphics[width=0.5\linewidth]{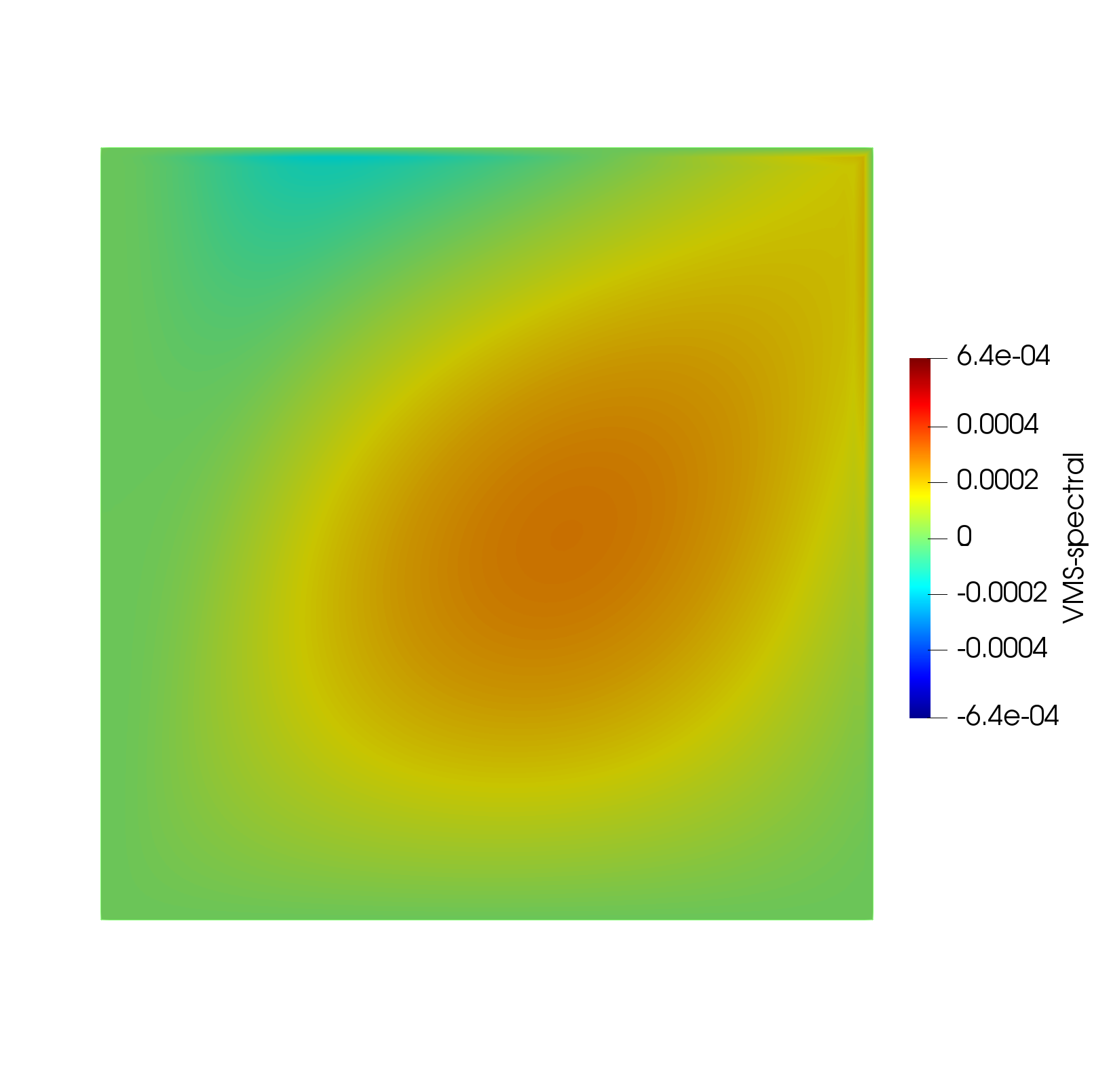}&\includegraphics[width=0.5\linewidth]{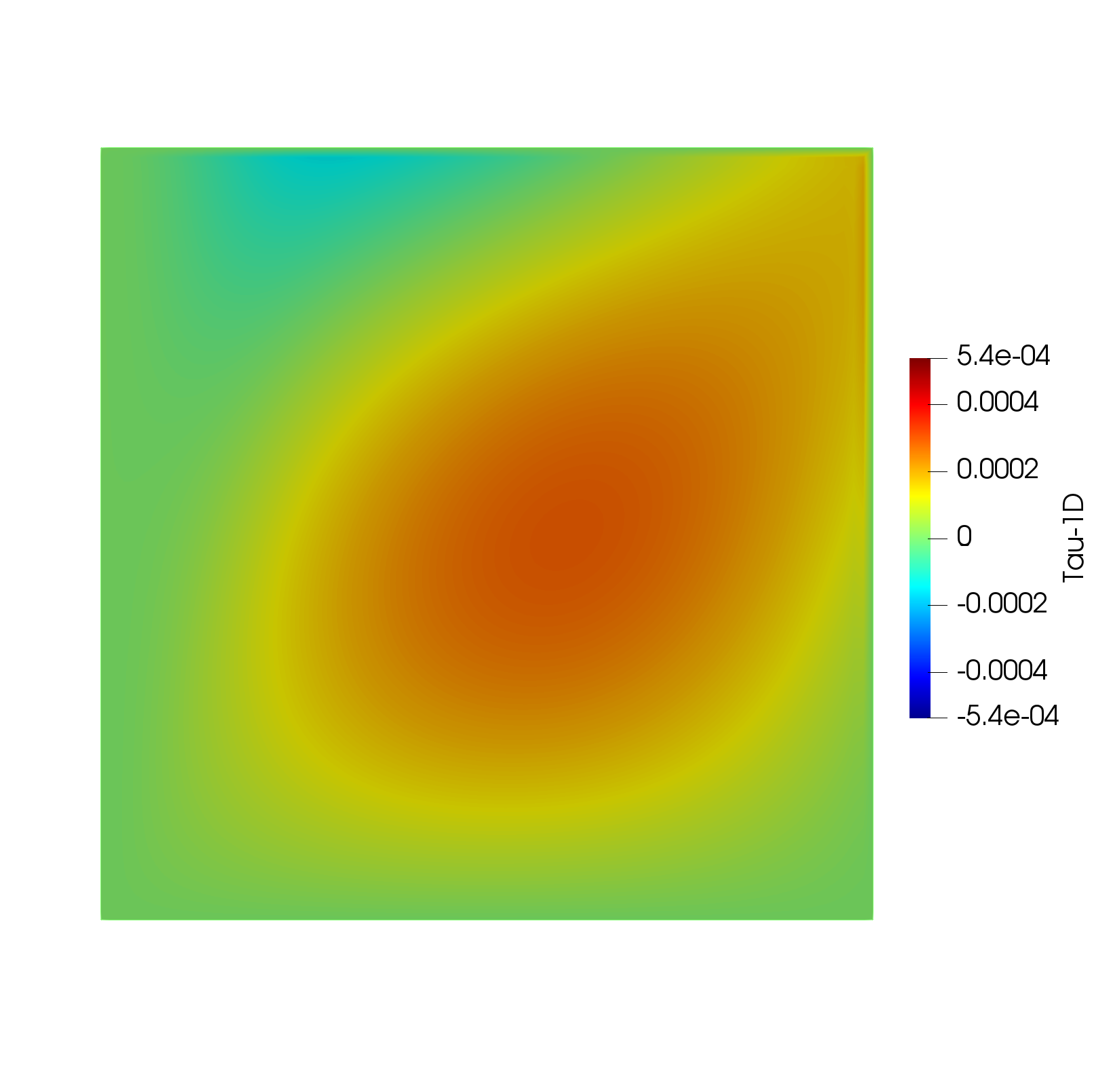}
\end{tabular}
\end{center}
\caption{\label{fig_taucte2}Representation of the solution obtained with the VMS-spectral method, panel (a) and using the stabilization coefficient $\tau_G$, panel (b) in the case $n=4.$}
\end{figure}

\begin{figure}[h!]
\begin{center}
\begin{tabular}{ll}(a)  Codina stabilized coefficients&(b)  Exact solution\\
\includegraphics[width=0.5\linewidth]{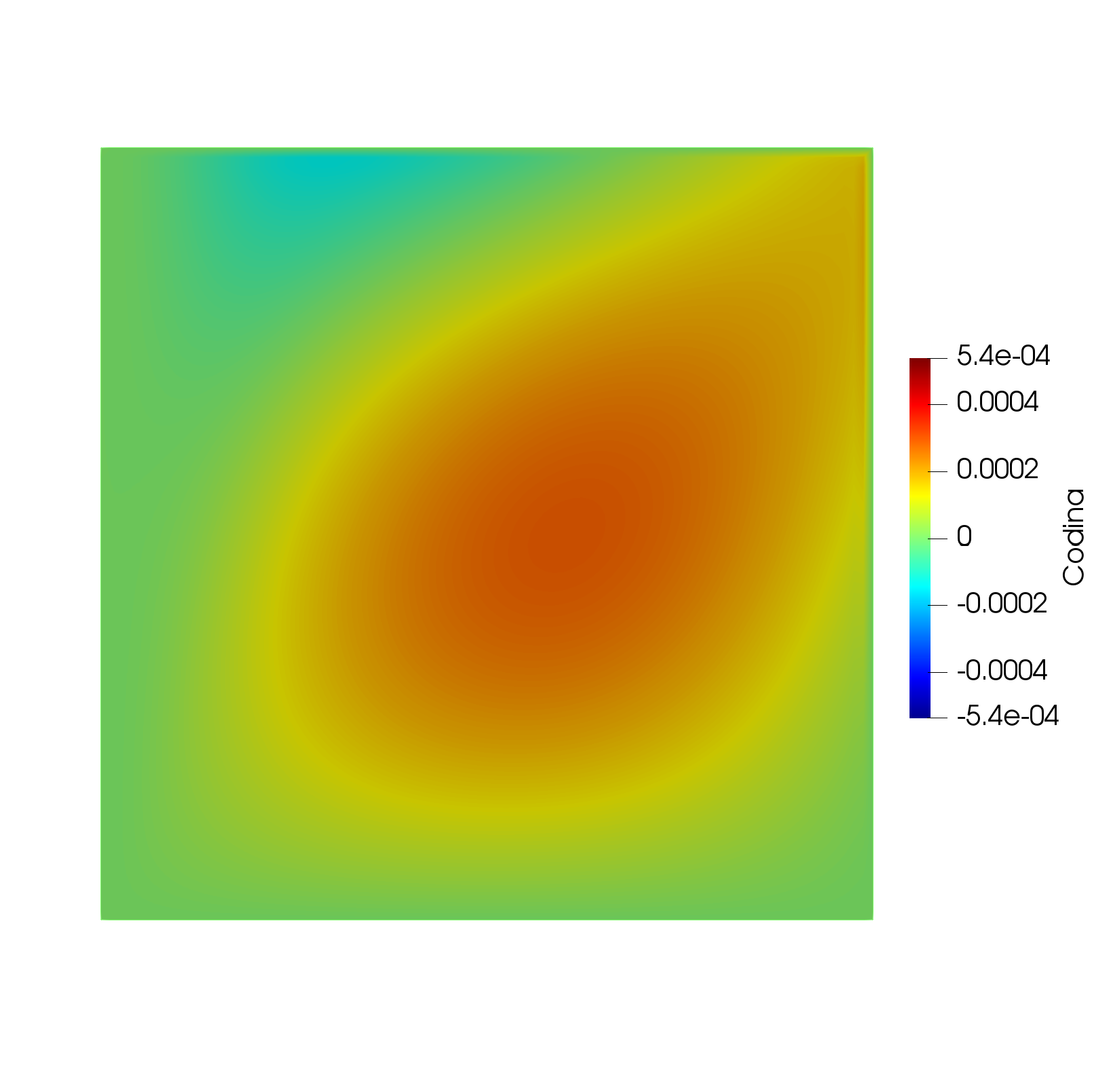}&\includegraphics[width=0.5\linewidth]{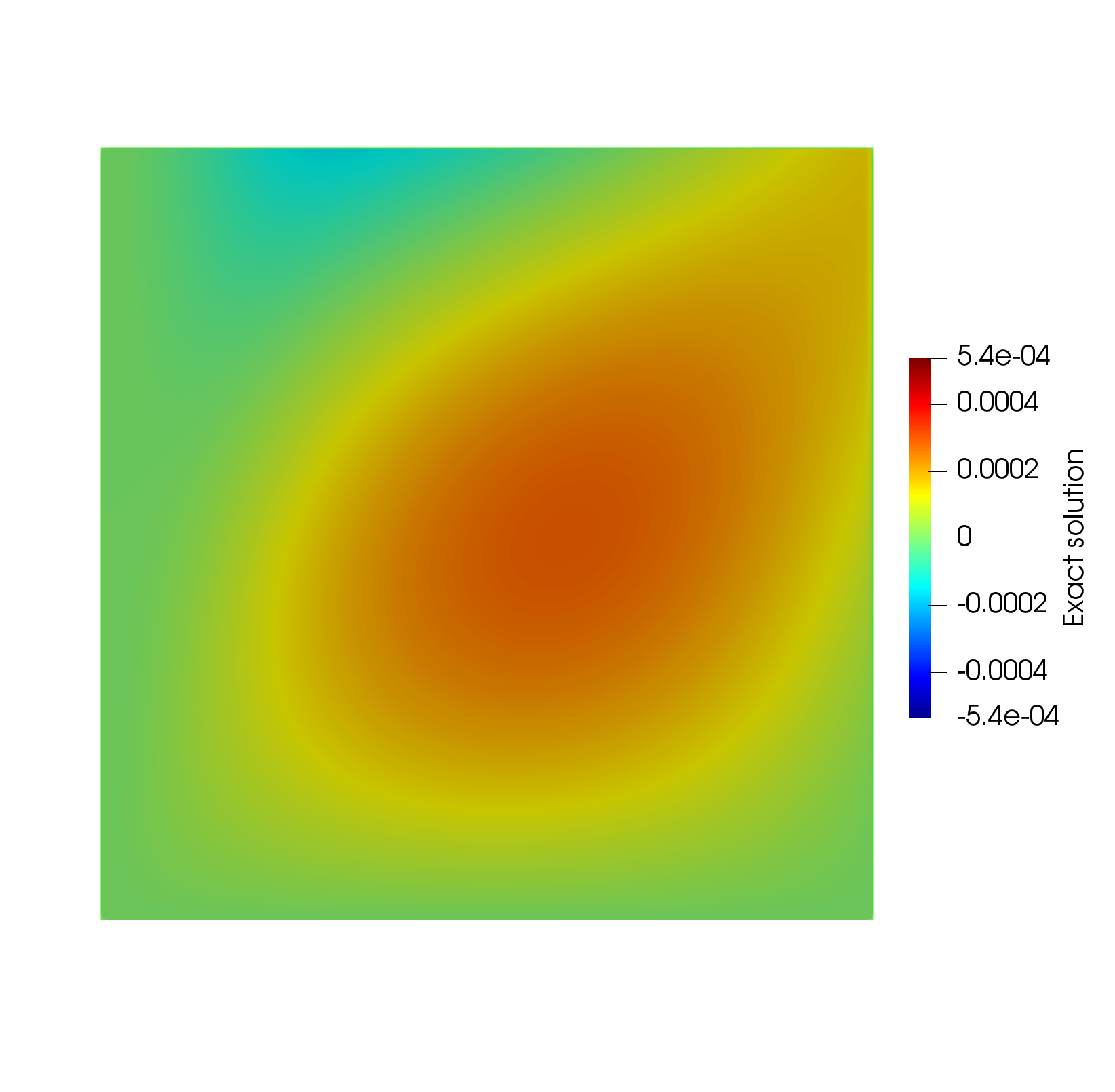}
\end{tabular}
\end{center}
\caption{\label{fig_solcte2} Representation of the solution obtained with $\tau_C$, panel (a) and the exact solution, panel (b) in the case $n=4.$}
\end{figure}

In particular, in Fig. \ref{fig_taucte2}, we represent solution obtained with the VMS-spectral method, panel (a) and using the stabilization coefficient $\tau_G$, panel (b) in the case $n=4,$ which corresponds to ${Pe_h}_1=2.18 $ and ${Pe_h}_2=6.72$. 
In Fig. \ref{fig_solcte2},  we represent solution obtained with the stabilization coefficient $\tau_C$, panel (a) and the exact solution, panel (b) in the case $n=4,$ which corresponds to ${Pe_h}_1=2.18 $ and ${Pe_h}_2=6.72$. 

\begin{figure}[h!]
\begin{center}
\begin{tabular}{ll}(a)   VMS-spectral stabilized coefficients&(b)  Generalized 1D stabilized coefficients\\
\includegraphics[width=0.5\linewidth]{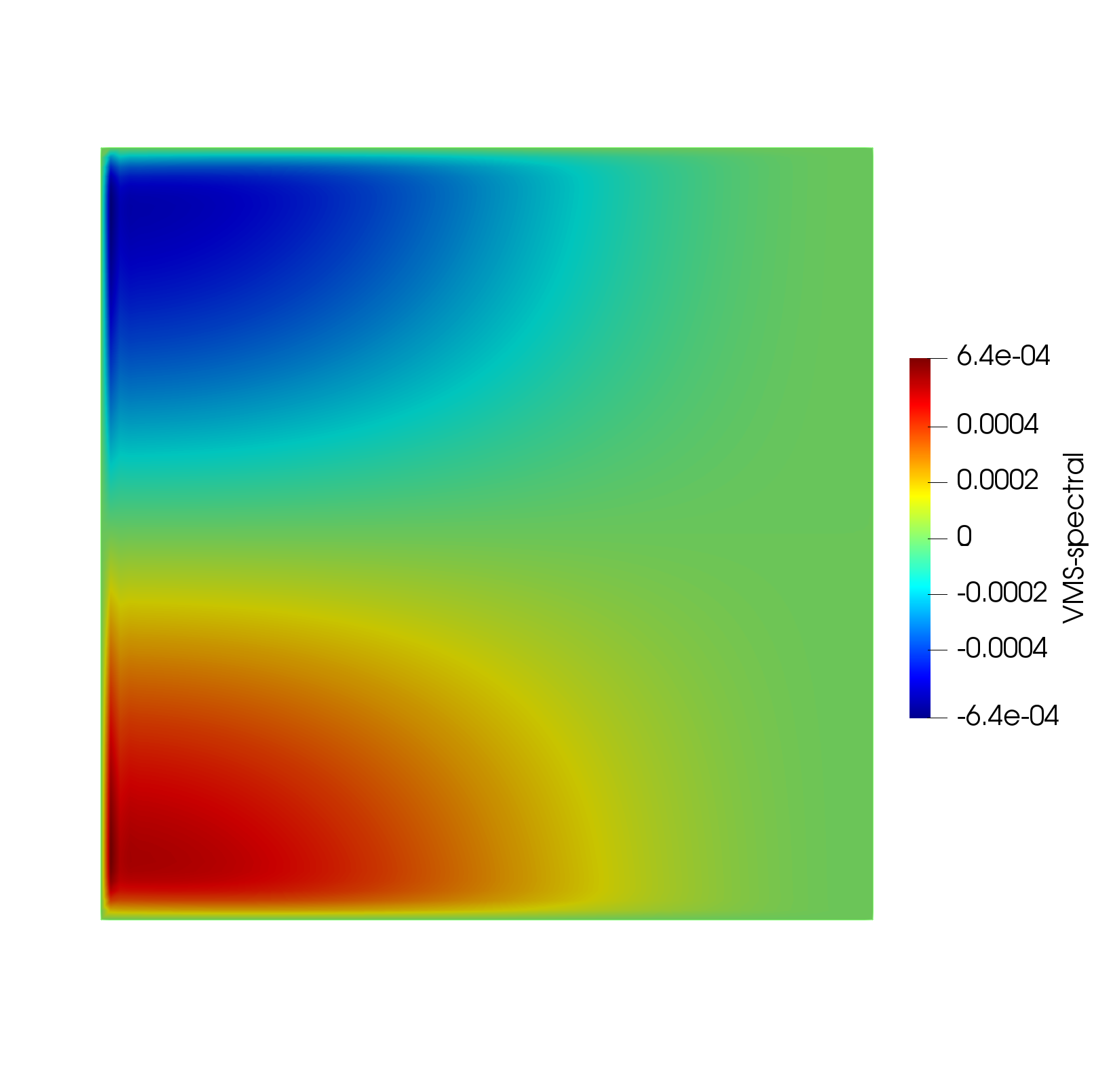}&\includegraphics[width=0.5\linewidth]{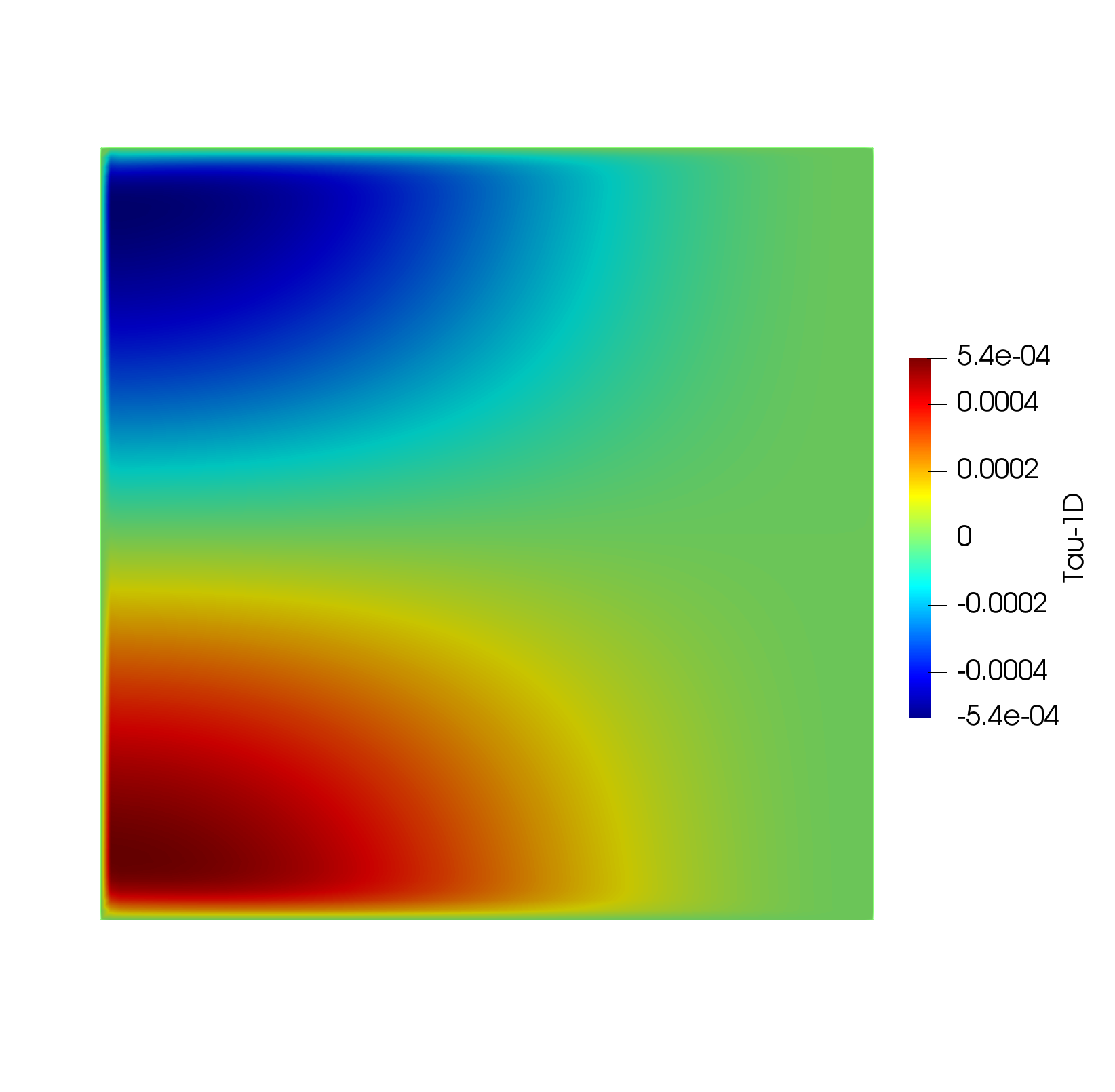}
\end{tabular}
\end{center}
\caption{\label{fig_taucte3}Representation of the solution obtained with the VMS-spectral stabilized coefficients, panel (a) and using the stabilization coefficient $\tau_G$, panel (b) in the case $n=10.$}
\end{figure}

\begin{figure}[h!]
\begin{center}
\begin{tabular}{ll}(a)  Codina stabilized coefficients&(b)  Exact solution\\
\includegraphics[width=0.5\linewidth]{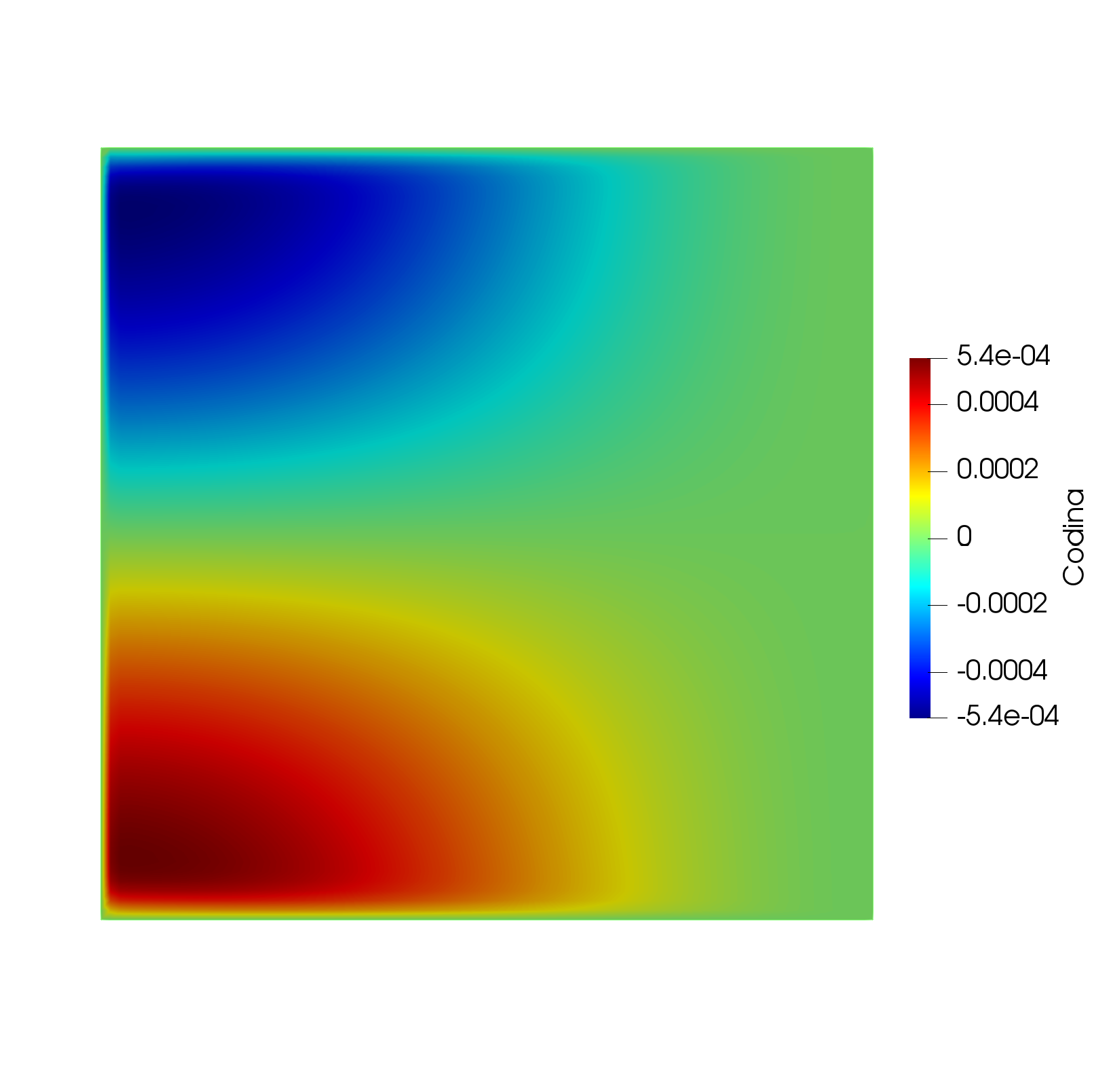}&\includegraphics[width=0.5\linewidth]{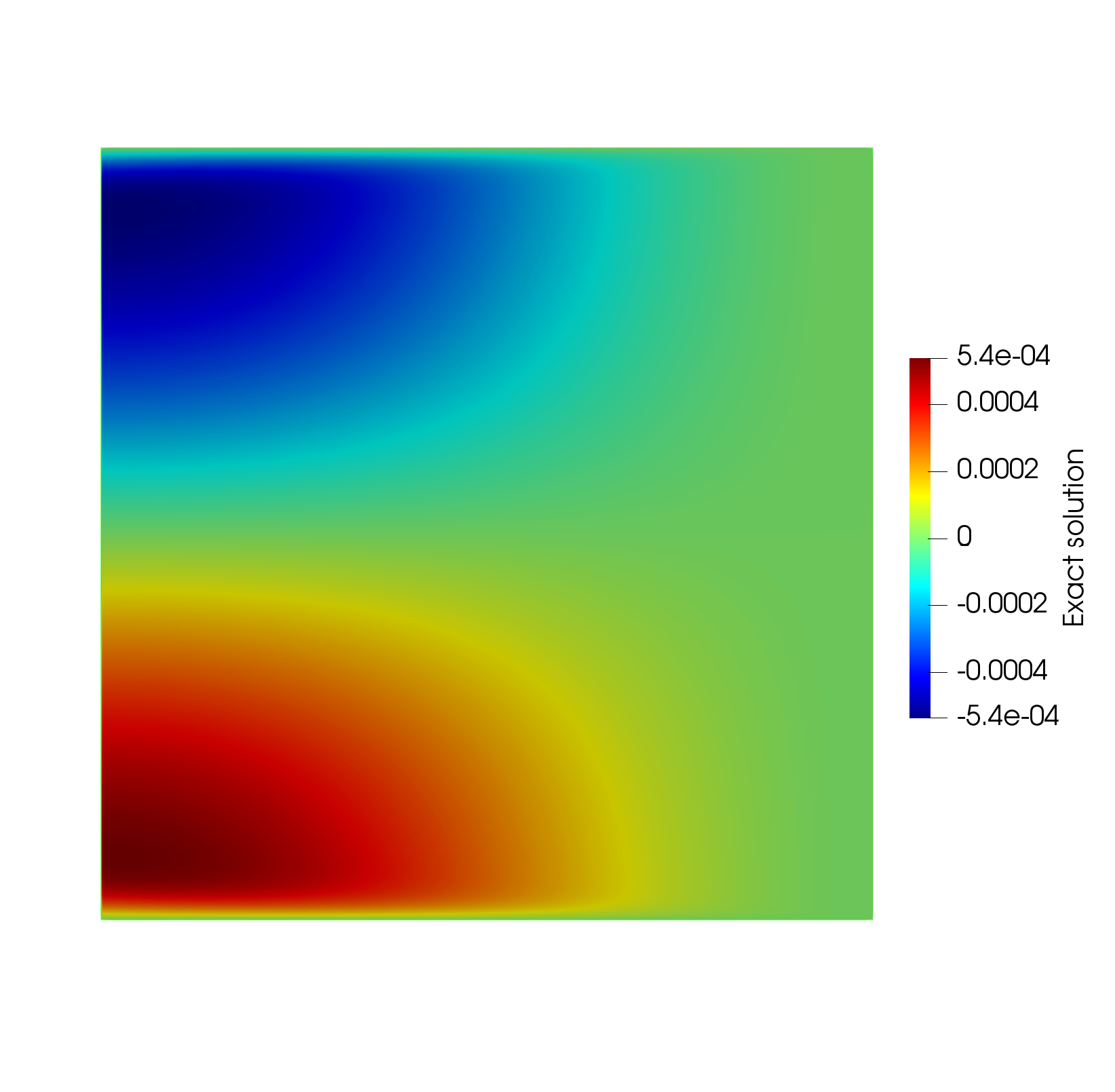}
\end{tabular}
\end{center}
\caption{\label{fig_solcte3} Representation of the solution obtained with $\tau_C$, panel (a) and the exact solution, panel (b) in the case $n=10.$}
\end{figure}

In Fig. \ref{fig_taucte3}, we represent solution obtained with the VMS-spectral stabilized coefficients, panel (a) and using the stabilization coefficient $\tau_G$, panel (b) in the case $n=10,$ which corresponds to ${Pe_h}_1=-7.07 $ and ${Pe_h}_2=0$. 
In Fig. \ref{fig_solcte3},  we represent solution obtained with the stabilization coefficient $\tau_C$, panel (a) and the exact solution, panel (b) in the case $n=10,$ which corresponds to ${Pe_h}_1=-7.07$ and ${Pe_h}_2=0$. 

\subsection{Tests with anisotropic velocity}
Here we consider a second test with anisotropic velocity, where we compare the VMS-spectral method with Codina stabilized coefficients.

\begin{figure}[h!]
\begin{center}
\includegraphics[width=0.6\linewidth]{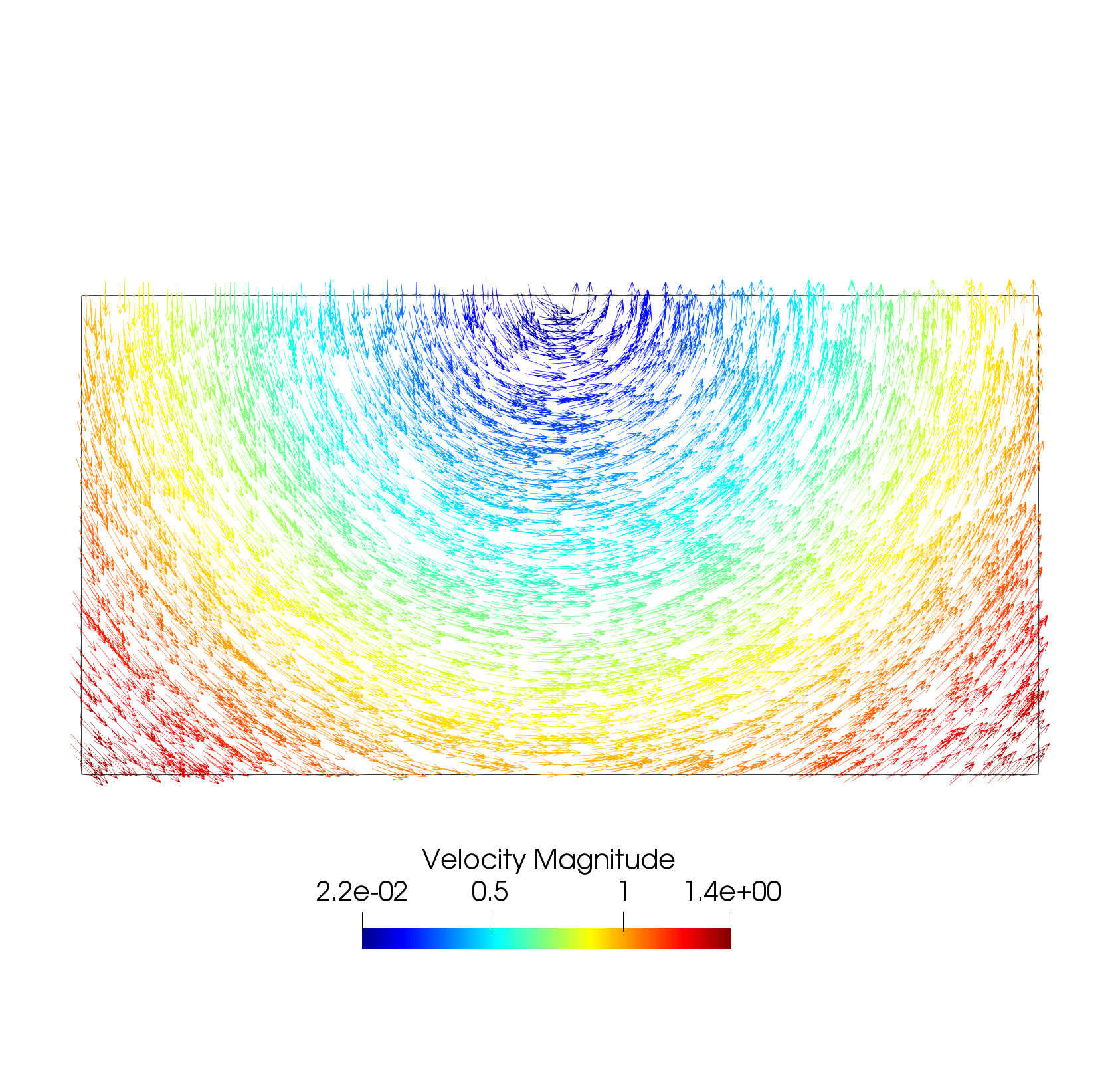}
\end{center}
\caption{\label{fig_vel} Representation of the velocity $\mathbf{a}(x,y)$.}
\end{figure}

Consider the advection-diffusion problem (\ref{EADS}) in the rectangle $[0,1]\times[0,1/2]$ with Dirichlet boundary conditions
with diffusion coefficient $\mu=10^{-3},$ source term $f(x,y)=1$ and velocity function
\begin{equation}\label{vel}
\mathbf{a}(x,y)=(a_1(x,y),a_2(x,y)),
\end{equation}
being
\begin{equation}\label{a1}
a_1(x,y)=\left\{\begin{array}{ll}
-0.1(y-0.5)&\mbox{if }\sqrt{x^2+y^2}<0.01, \\ \noalign{\smallskip}
-2(y-0.5)&\mbox{if }\sqrt{x^2+y^2}\geq 0.01,
\end{array}\right.
\end{equation}
and
\begin{equation}\label{a2}
a_2(x,y)=\left\{\begin{array}{ll}
0.1(x-0.5)&\mbox{if }\sqrt{x^2+y^2}<0.01, \\ \noalign{\smallskip}
2(x-0.5)&\mbox{if }\sqrt{x^2+y^2}\geq 0.01.
\end{array}\right.
\end{equation}
In Fig. \ref{fig_vel}, we represent velocity vector field $\mathbf{a}$ given in (\ref{vel})-(\ref{a2}). 

In Table \ref{tab:table1}, we provide for each $N=1/h,$ being $h$ the diameter of the grid elements: the maximum P\'eclet number  ${Pe_h}$, the ranges of values of ${Pe_h}_1$ and ${Pe_h}_2$ (note that these P\'eclet numbers depend on the component of the velocity and can have a negative sign),
the error in $L^2$ and $L^\infty$ norms of the stabilized solution through the VMS-spectral stabilized coefficients (with subindex $A$) and through the Codina stabilized coefficients (with subindex $C$),
  \begin{center}
       \label{tab:table1}
   \begin{tabular}{c c c c c c c c}  
   \hline\hline                        
   N&${Pe_h}$-max&${Pe_h}_1-$range &${Pe_h}_2-$range  &$L^2_A$ & $L^2_C$ & $L^\infty_A$ & $L^\infty_C$ \\    [0.5ex]
   \hline                   
   50 &6.97&(0.03,4.96) &(-4.93,4.93) & 5.080e-03 & 9.972e-03 & 0.1317 & 0.1720  \\  
   100 &3.51&(8.33e-05,2.49) & (-2.48,2.48)& 1.717e-03 & 5.15e-03 & 0.065 & 0.1149\\ 
   200 &1.75&(-2.08e-05,1.24) & (-1.24,1.24)& 4.83e-04 & 2.18e-03  &  0.025 & 0.0604\\ 
   [0.5ex]       
   \hline     
   \end{tabular}
    \end{center}
    
\begin{figure}[h!]
\begin{center}
\begin{tabular}{ll}(a)  VMS-spectral stabilized coefficients&(b) Codina stabilized coefficients\\
\includegraphics[width=0.5\linewidth]{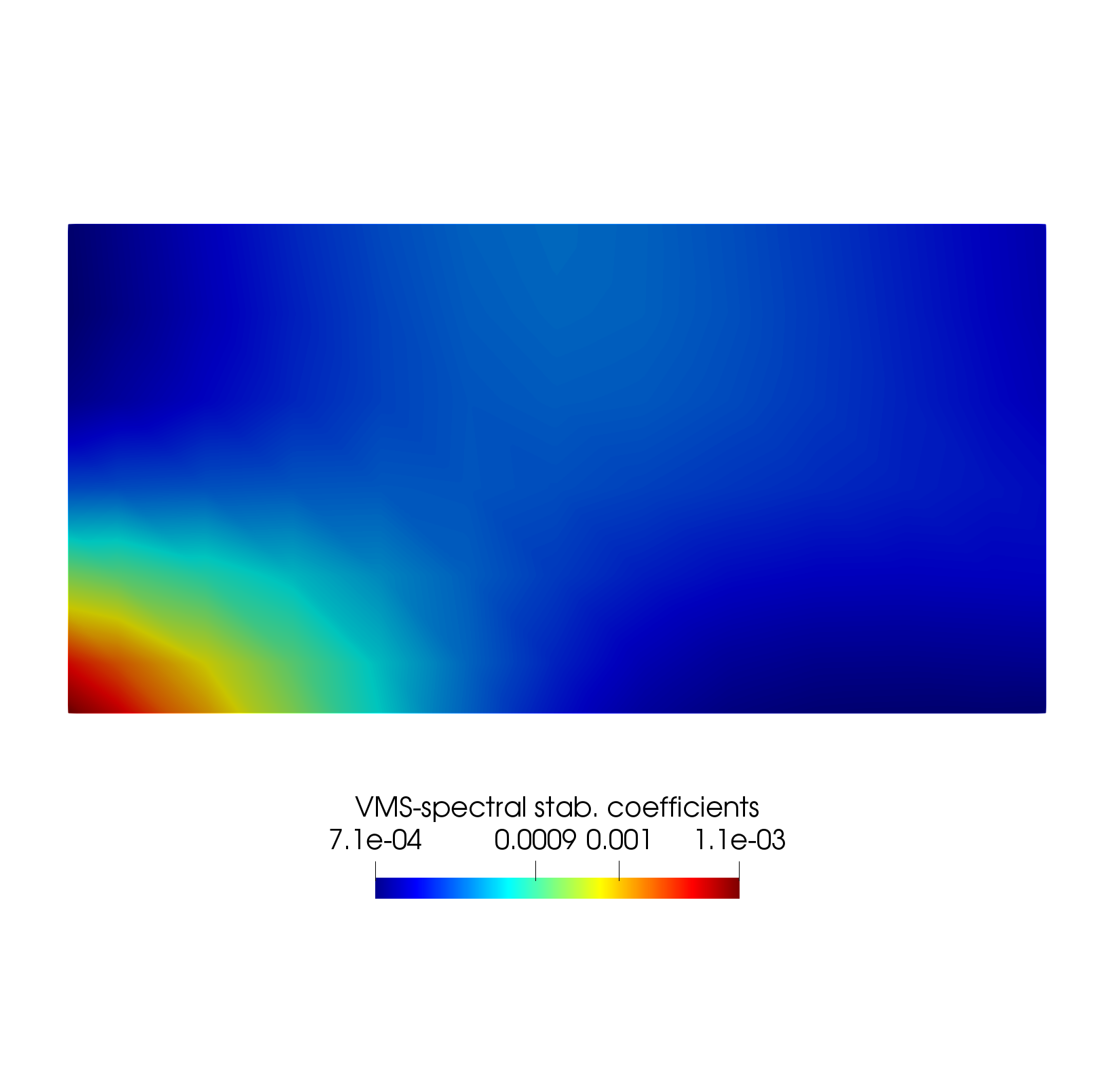}&\includegraphics[width=0.5\linewidth]{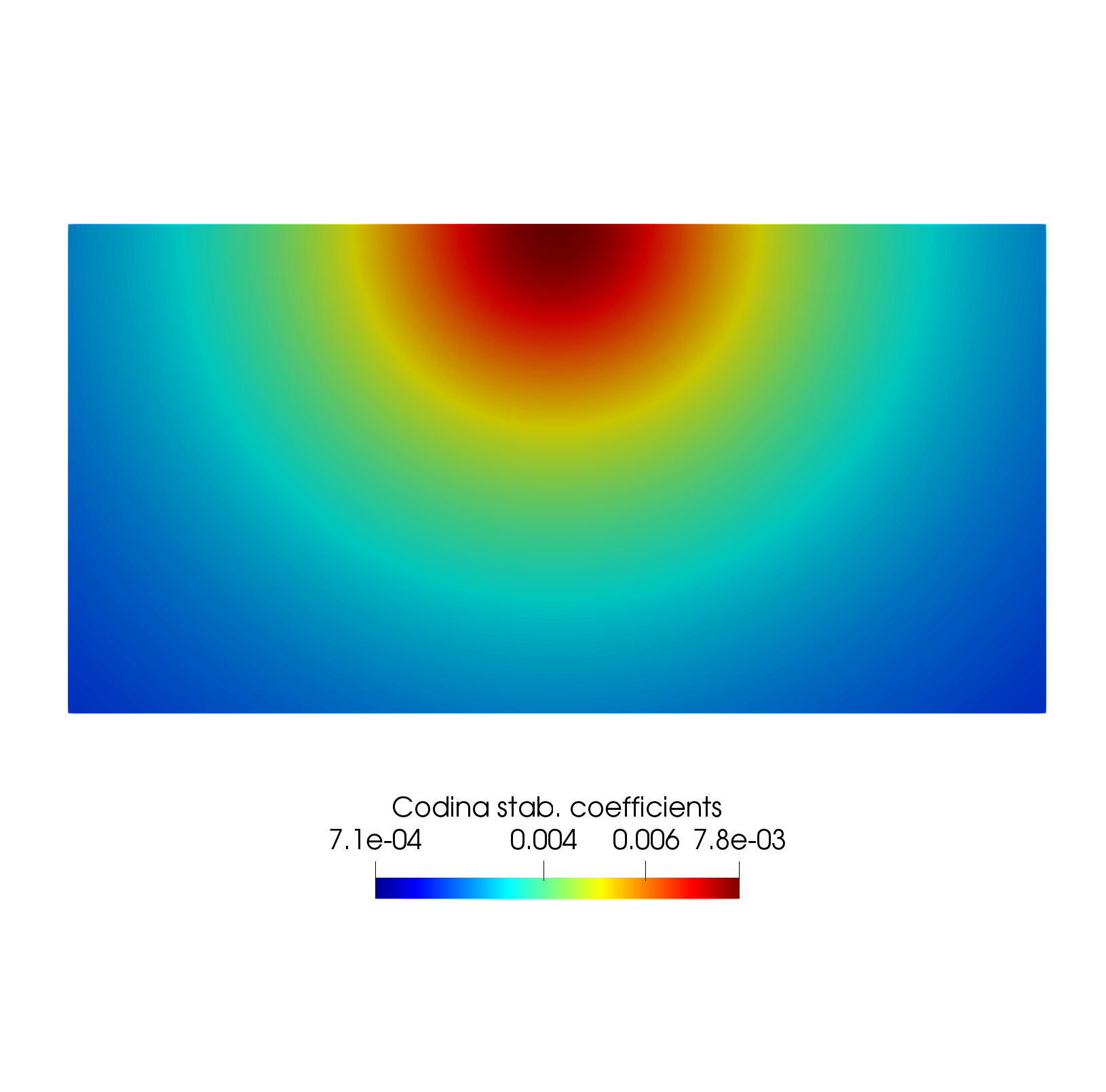}
\end{tabular}
\end{center}
\caption{\label{fig_taus}Representation of the stabilization coefficients obtained with: the VMS-spectral stabilized coefficients, (panel (a)) and the Codina stabilized coefficients, (panel (b)).}
\end{figure}

In Fig. \ref{fig_taus}, we represent in the case $N=100$ the stabilization coefficients obtained with the VMS-spectral stabilized coefficients, panel (a) and through Codina stabilized coefficients, panel (b). 

\begin{figure}[h!]
\begin{center}
\includegraphics[width=0.7\linewidth]{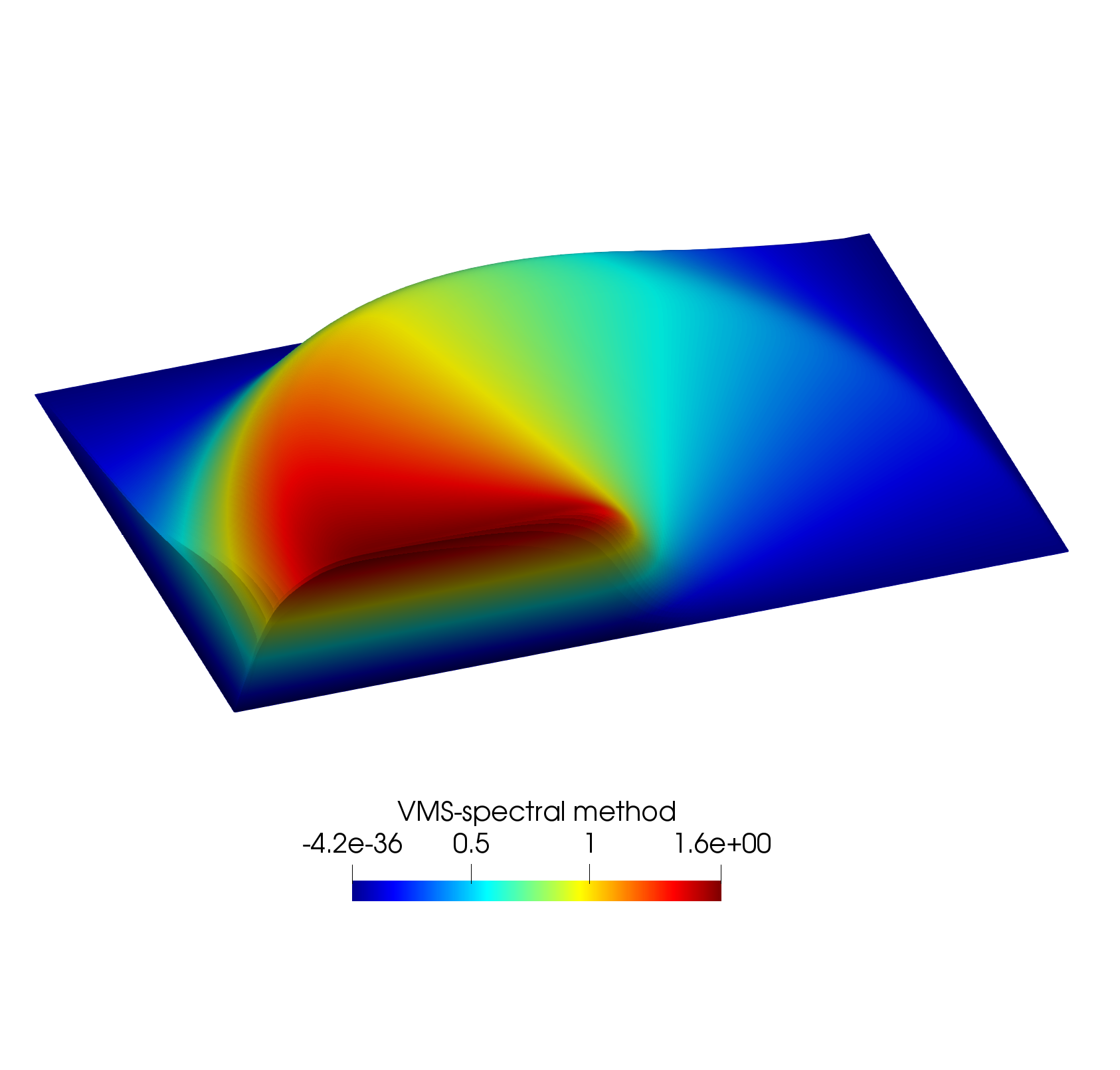} \\
\end{center}
\caption{\label{fig_sol} Representation of the solution obtained with the VMS-spectral method. }
\end{figure}

Finally, in Fig. \ref{fig_sol}, we represent in the case $N=100$ the solution obtained with the VMS-spectral stabilized coefficients.
%
 \subsection{Flow around a cylinder}
The tests considered in previous subsections have been performed in regular meshes of isosceles right triangles, for which the theory has been developed.  Finally, in this subsection, we check the reliability of the method in a more realistic case with an irregular mesh.

\begin{figure}[h!]
\begin{center}
\includegraphics[width=0.7\linewidth]{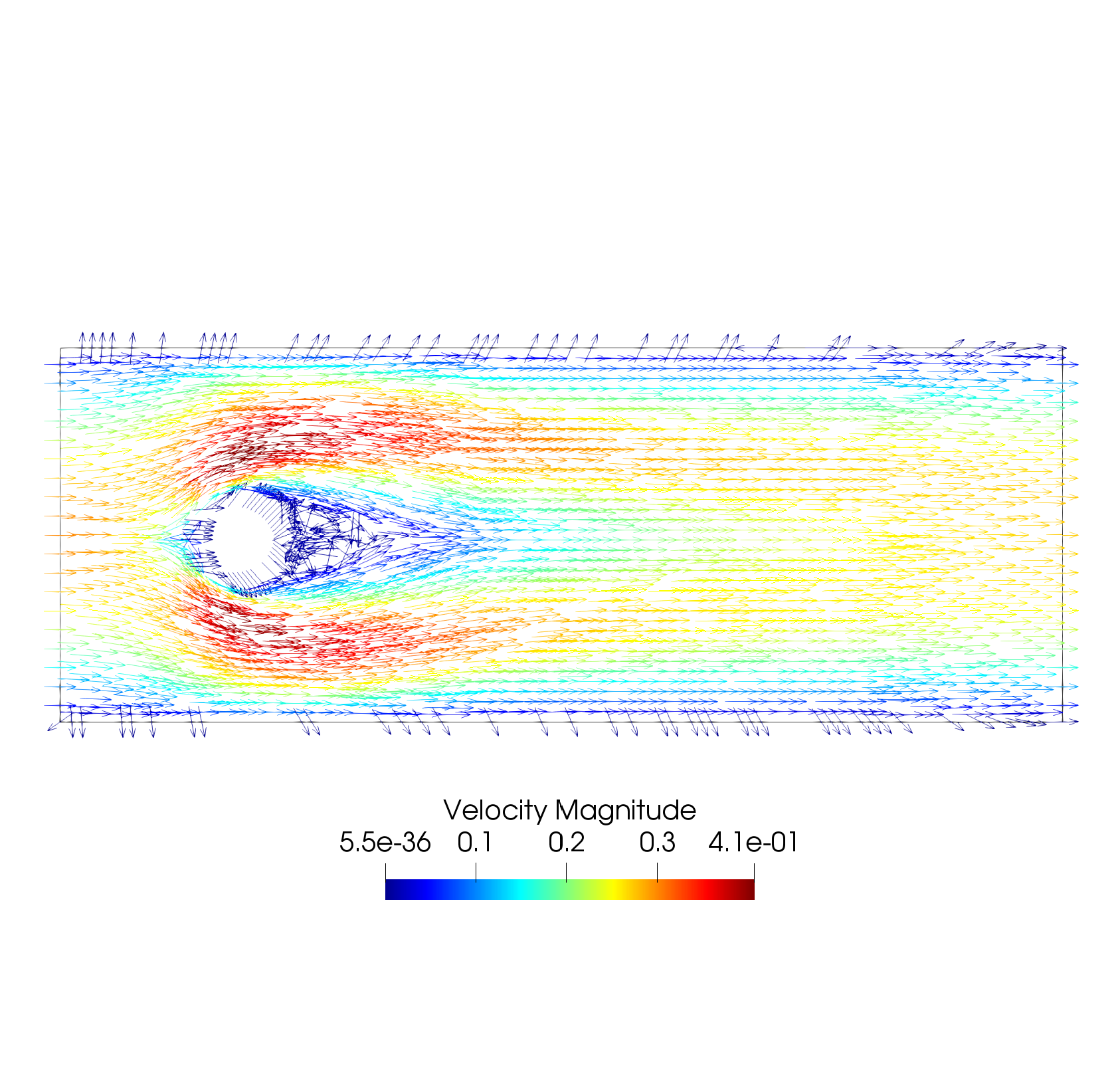}
\end{center}
\caption{\label{fig_vel_cil} Representation of the velocity $\mathbf{a}(x,y)$.}
\end{figure}

In this test, we have previously computed the stationary state of a fluid with Reynolds number 100 around a cylinder. We use this velocity as the transport velocity $\mathbf{a}(x,y)$ to solve the advection-diffusion problem  (\ref{EADS})  for a passive agent.  
In Fig. \ref{fig_vel_cil}, we represent the velocity vector field $\mathbf{a}$.

Moreover, we consider diffusion coefficient $\mu=10^{-3}$ and source term $f(x,y)=1$.  
The maximum P\'eclet number is ${Pe_h}=2.1,$ while the directional P\'eclet numbers ${Pe_h}_1\in(0,4.08)$ and
${Pe_h}_2\in(-1,1.01).$  

\begin{figure}[h!]
\begin{center}
\begin{tabular}{ll}(a)  VMS-spectral stabilized coefficients&(b) Codina stabilized coefficients\\
\includegraphics[width=0.5\linewidth]{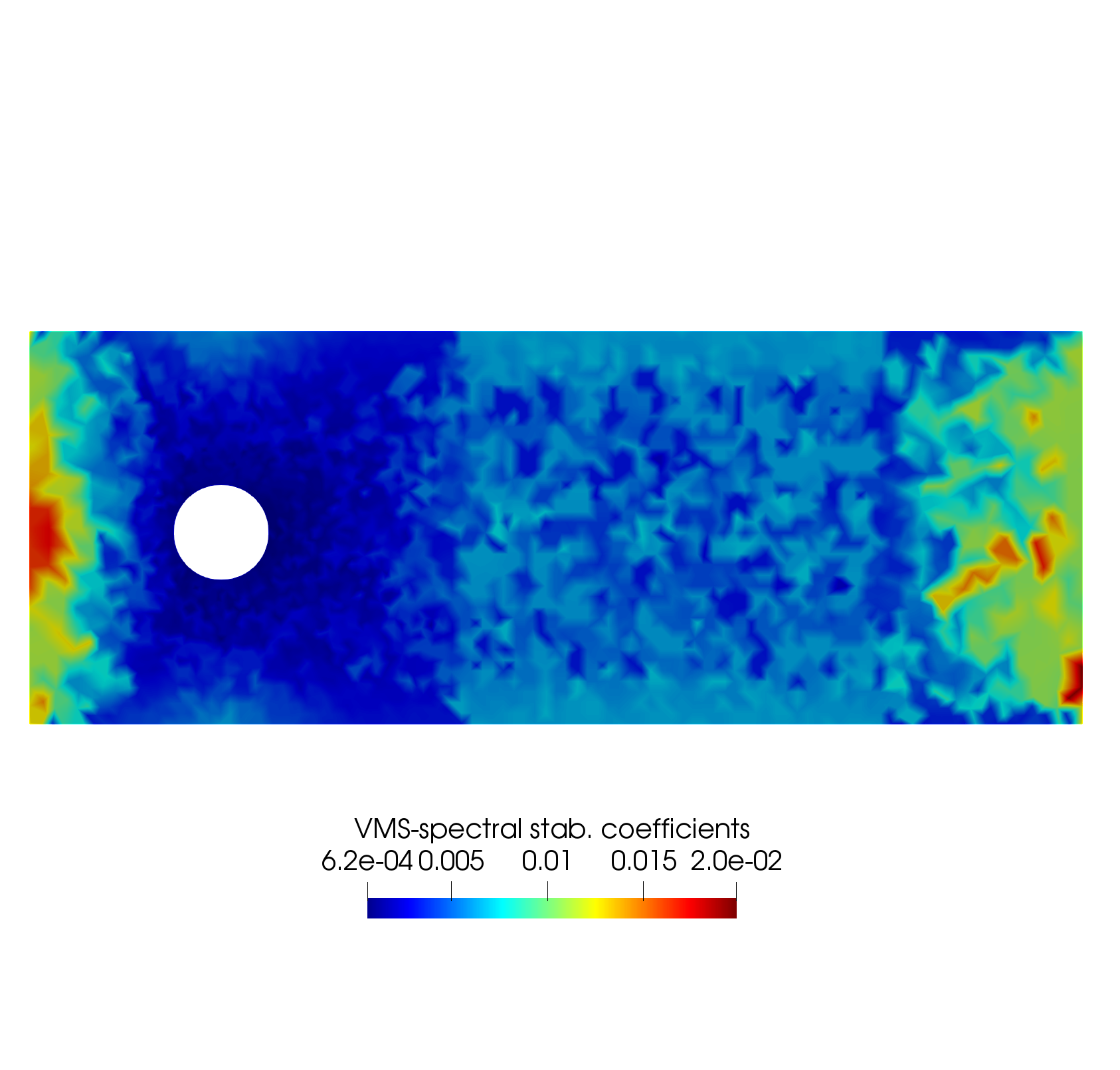}&\includegraphics[width=0.5\linewidth]{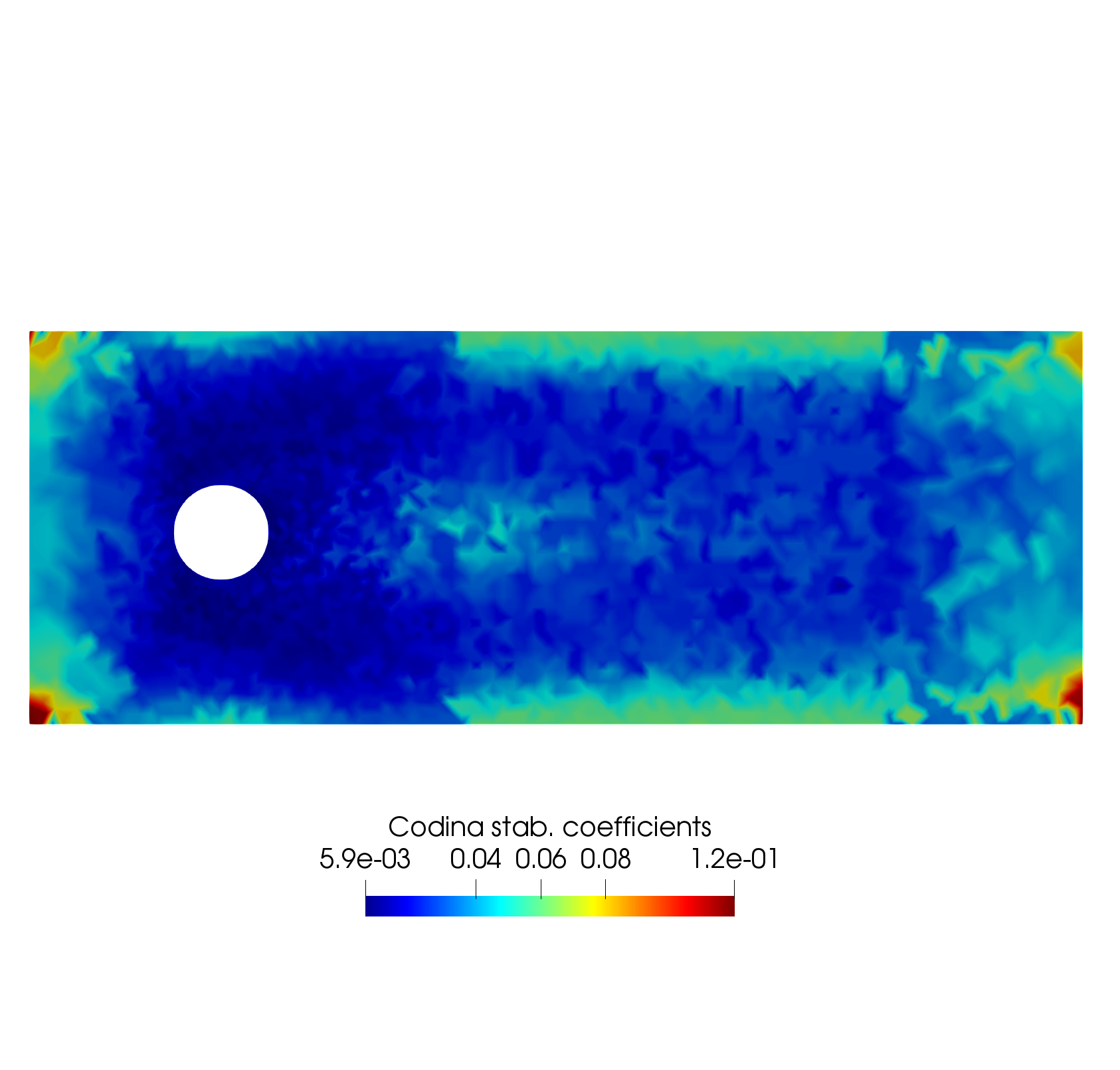}
\end{tabular}
\end{center}
\caption{\label{fig_taus_cil}Representation of the stabilization coefficients obtained with: the VMS-spectral stabilized coefficients, (panel (a)) and the Codina stabilized coefficients, (panel (b)). }
\end{figure}

In Fig. \ref{fig_taus_cil}, we represent the stabilization coefficients obtained with the VMS-spectral stabilized coefficients, (panel (a)) and through Codina stabilized coefficients, (panel (b)). 

\begin{figure}[h!]
\begin{center}
\includegraphics[width=0.7\linewidth]{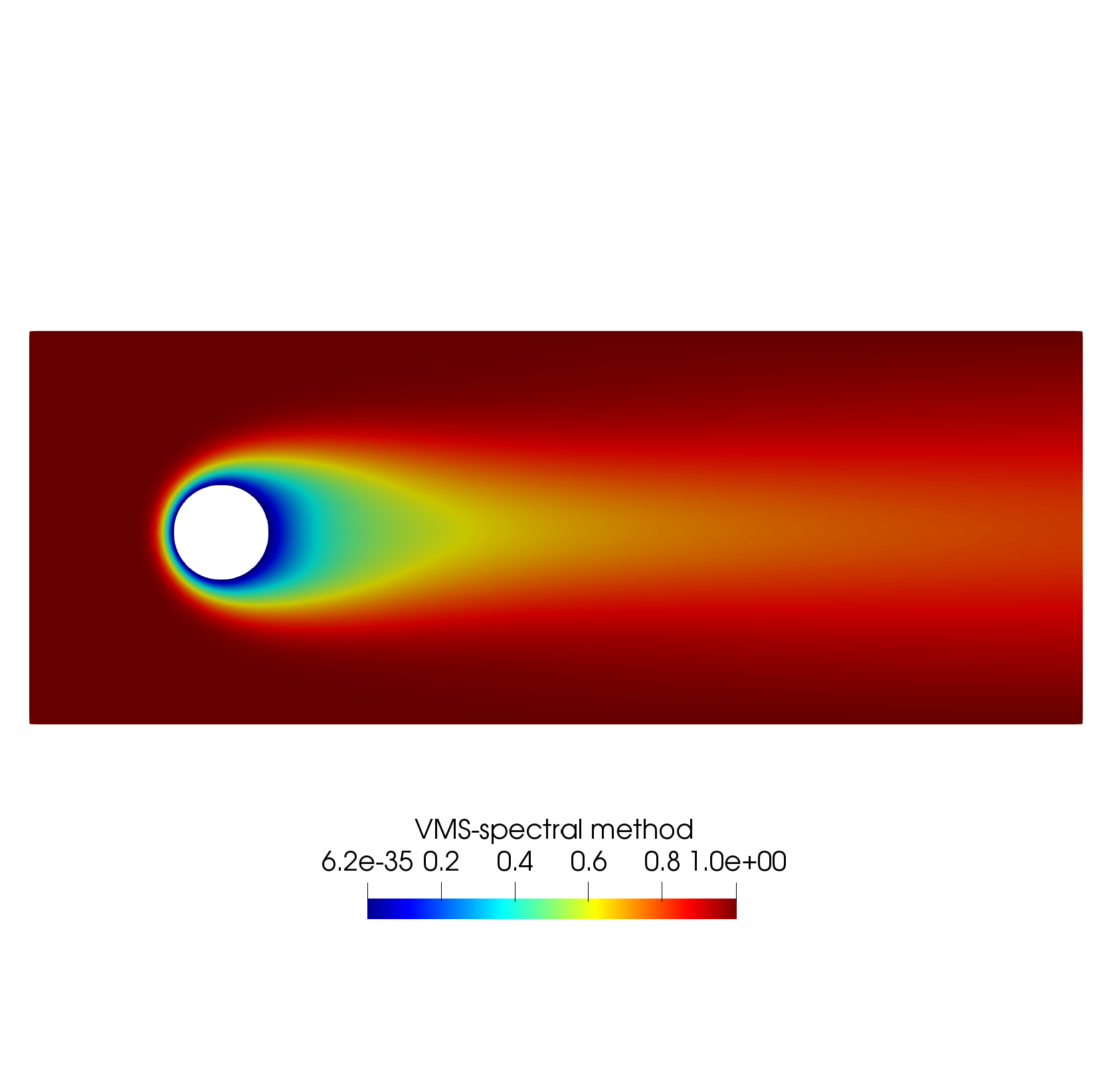}
\end{center}
\caption{\label{fig_sol_cil} Representation of the solution obtained with the VMS-spectral method.}
\end{figure}

In Fig. \ref{fig_sol_cil}, we represent the solution obtained with the VMS-spectral method.

\begin{figure}[h!]
\begin{center}
\begin{tabular}{ll}(a) VMS-spectral stabilized coefficients&(b)  Codina stabilized coefficients\\
\includegraphics[width=0.5\linewidth]{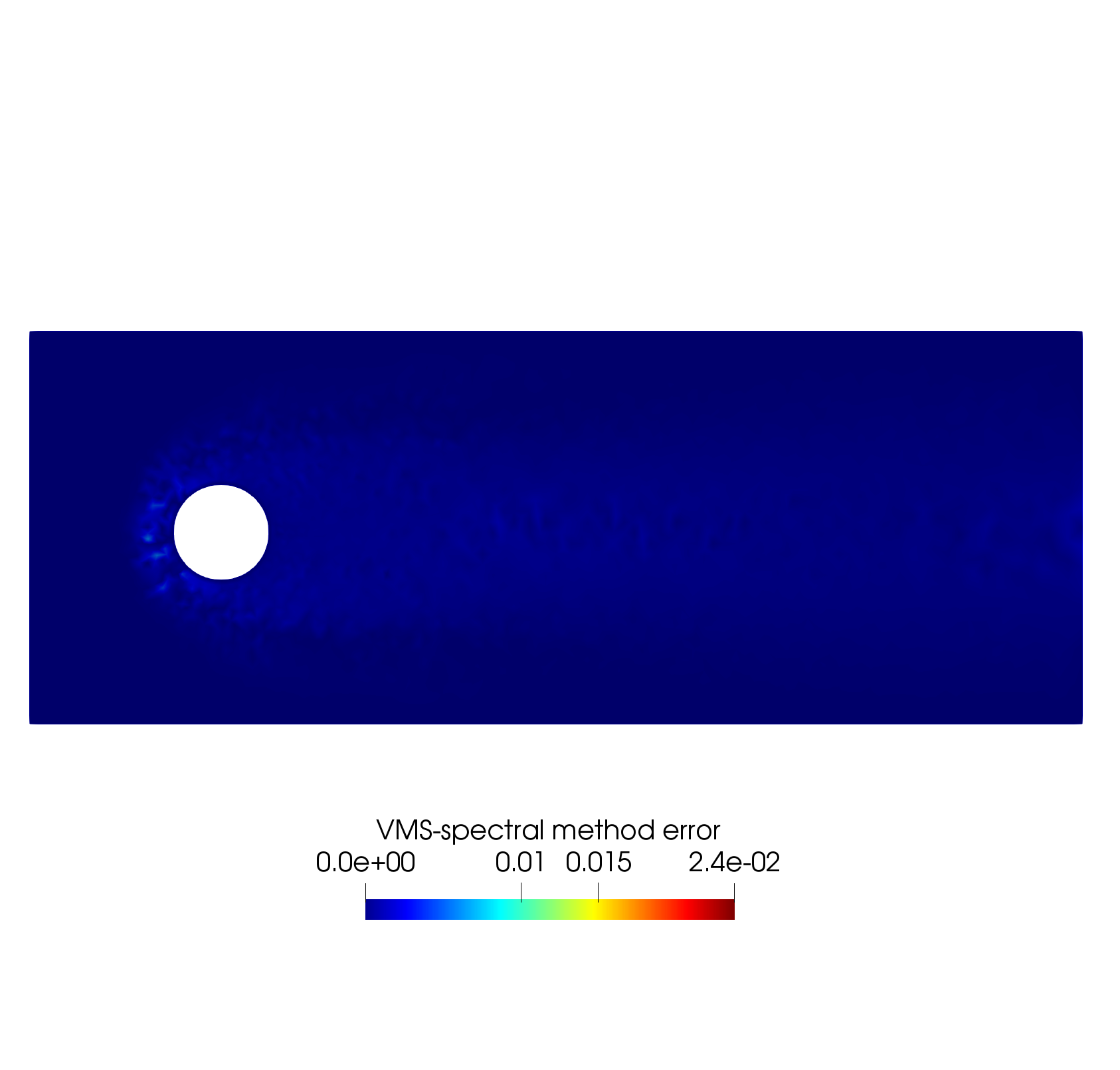}&\includegraphics[width=0.5\linewidth]{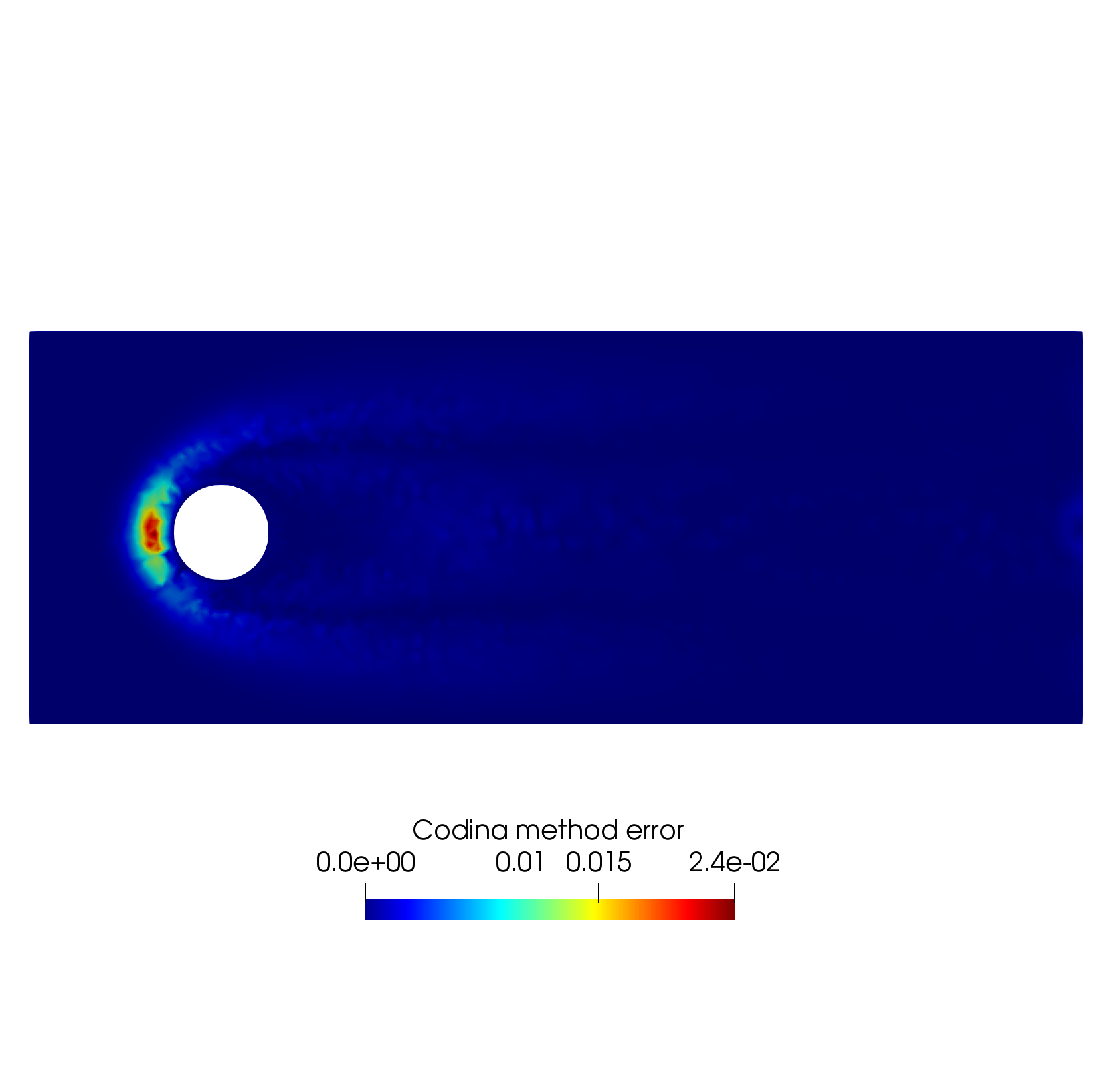}
\end{tabular}
\end{center}
\caption{\label{fig_error_cil} Representation of the errors obtained with the VMS-spectral stabilized coefficients, (panel (a)) and with Codina stabilized coefficients, (panel (b)). }
\end{figure}

In Fig. \ref{fig_error_cil}, we represent the errors obtained with the VMS-spectral stabilized coefficients, (panel (a)) and with Codina stabilized coefficients, (panel (b)). 
We have obtained the errors 5.35e-04  and 1.2e-03 in $L^2-$norm with the VMS-spectral stabilized coefficients and Codina stabilized coefficients, respectively. Regarding the $L^\infty-$norm, we have obtained
an error of 2.39e-02 using the VMS-spectral stabilized coefficients and of 3.81e-02 with Codina stabilized coefficients.

\section{Conclusions and perspectives}\label{sec:conclusions}

In this paper, we have extended to two dimensional problems the VMS-spectral method developed in \cite{ChaconDia} for 1D elliptic problems, and we have applied it to the advection-diffusion problem.    

For piecewise affine finite element discretizations, the Spectral VMS method can be cast as a standard VMS method with approximate stabilized coefficients. As the eigen-pairs of the spectral sub-grid problems can be explicitly computed, analytical expressions for the stabilized coefficients have been obtained. In particular, for 2D advection-diffusion equations, these can be computed in an off-line phase, and then interpolated in the on-line phase. 

In order to test the reliability of the method, we have provided numerical tests with constant and anisotropic velocity with intermediate P\'eclet number values in regular and irregular meshes. We have compared the results provided by the spectral method, to those obtained with generalized 1D stabilization coefficient \cite{christie,johnkno} and with the stabilization coefficients obtained by Codina \cite{codina12}.

Regarding future work, let us note that in the 1D case, apart from the advection-diffusion problem, the VMS-spectral method has been applied to the advection-diffusion-reaction problem with Dirichlet boundary conditions \cite{ChaconFernandez}. This has  allowed us
to extend the method to parabolic problems, with application to evolution advection-diffusion equations \cite{ChaconFernandez2}. As in the 2D case, the eigen-pairs of the 2D advection-diffusion-reaction operator with Dirichlet boundary conditions in structured grids can be also computed from those of the Laplace operator, then the VMS method can be directly applied to 2D advection-diffusion-reaction problems. Therefore it could be adapted to two-dimensional evolutive problems. 

Finally, as another perspective, as we have already remarked before, while in 1D case the decomposition of the Hilbert space into the large scales space plus the small scales space is exact, this is not the case in 2D. Therefore,
in order to have an exact decomposition of the 2D space, subscales on the element boundaries need to be taken into account.

\section*{Acknowledgements}
This research is partially supported by Ministerio de Econom\'ia
y Competitividad under grant MTM2015-64577-C2-1-R. S.
Fern\'andez-Garc\'ia is supported by the University of Seville
VPPI-US.

\end{document}